\documentclass[12pt]{paper}
\usepackage{graphicx}
\usepackage{amsmath}
\usepackage{url}
\usepackage{verbatim}
\usepackage{epsf}
\usepackage{latexsym}
\usepackage{array,tabu,multirow,multicol,longtable}
\usepackage{mathrsfs}
\usepackage{authblk}
\usepackage{amsfonts}
\usepackage{amssymb}
\usepackage{moreverb,footmisc,amsthm}%,fnpct}
\usepackage{color}
\usepackage[left=2.5cm,right=2cm,top=2cm,bottom=2cm]{geometry}
\usepackage{hyperref}

\hypersetup{
	colorlinks=true,
	linkcolor=blue,
	urlcolor=magenta,
	citecolor=red,
}
\newcommand*{\QED}{\hfill\ensuremath{\square}}%

\newenvironment{Proof}{\vspace{1.5ex}{\sc Proof}. }{\vspace{2ex} $\QED$}

\newtheorem{Theorem}{Theorem}
\newtheorem{Corollary}{Corollary}
\newtheorem{Definition}{Definition}

\newtheorem{Lemma}{Lemma}
\newtheorem{Proposition}{Proposition}
\newtheorem{Remark}{Remark}

\newtheorem{Assumption}{Assumption}
\def\Nb{\mathbb{N}}

\def\Rb{\mathbb{R}}

\def\vg0{\mathbf{0}}
\def\to{\rightarrow}
\def\To{\longrightarrow}

\def\F{\mathscr{F}}

\def\P{\mathscr{P}}

\def\norm#1{
    \left\|
    #1\right\|
    }

\def\tr#1{\mathop{\mathrm{tr}}#1}
\def\span#1{\mathop{\mathrm{span}}#1}

\def\bigO#1{\mathcal{O}\left(#1\right)}
\def\smallo#1{o\left(#1\right)}
\def\argmin{\mathop{\mathrm{argmin}}}

\def\ud{\,\mathrm{d}}
\def\ran{\mathop{\mathrm{Ran}}}

\begin{document}

\raggedbottom 
\addtolength{\topskip}{0pt plus 10pt}

%-------------------
\title{{\bf \begin{center}
A new regularization method for linear exponentially ill-posed problems
\end{center}}}
\author{Walter Cedric {\sc Simo Tao Lee}\thanks{Institut de Math\'ematiques de Toulouse, Universit\'e Paul Sabatier, Toulouse, France \\
Email: wsimotao@math.univ-toulouse.fr}}
%\date{September 23, 2019}
%-------------------
\maketitle

\abstract 
This paper provides a new regularization method which is particularly suitable for linear exponentially ill-posed problems. Under logarithmic source conditions (which have a natural interpretation in terms of Sobolev spaces in the aforementioned context), concepts of qualifications as well as order optimal rates of convergence are presented. Optimality results under general source conditions expressed in terms of index functions are also studied. 
 
 Finally, numerical experiments on three test problems attest the better performance of the new method compared to the well known Tikhonov method in instances of exponentially ill-posed problems.
 \newline


\section{Introduction}
In this paper, we are interested in the solution to the equation
\begin{equation}
\label{main equation}
 T x = y
\end{equation}
where $T : X \to Y$ is a linear bounded operator between two infinite dimensional Hilbert spaces~$X$ and~$Y$ with non-closed range. The data~$y$ belongs to the range of~$T$ and we assume that we only have approximated data~$y^\delta$ satisfying 
\begin{equation}
\label{noise level equation}
||y^\delta - y || \leq \delta.
\end{equation}
   In such a setting, Equation \eqref{main equation} is ill-posed in the sense that the Moore Penrose generalized inverse~$T^\dagger$ of $T$ which maps~$y$ to the best-approximate solution~$x^\dagger$ of \eqref{main equation} is not continuous. Consequently a little perturbation on the data~$y$ may induce an arbitrarily large error in the solution~$x^\dagger$. 
Instances of such ill-posed inverse problems are encountered in several fields in applied sciences among which: signal and image processing, computer tomography, immunology, satellite gradiometry, heat conduction problems, inverse scattering problems, statistics and econometrics to name just a few (see, e.g.  \cite{book_engl,book_groetsch,reg_exp_ipp,book_kirsch,book_murio}).
   As a result of the ill-posedness of Equation \eqref{main equation}, a regularization method needs to be applied in order to recover from the noisy data~$y^\delta$ a stable approximation~$x^\delta$ of the solution~$x^\dagger$. A regularization method can be regarded as a family of continuous operators $R_\alpha : Y \to X$ such that there exists a function 
$\Lambda : \Rb_+ \times Y \to \Rb_+ $ satisfying the following: for every~$y \in \mathcal{D}(T^\dagger) \subset Y $ and~$y^\delta \in Y$ satisfying \eqref{noise level equation}
\begin{equation}
\label{convergence reg method}
 R_{\Lambda(\delta,y^\delta)} y^\delta \to x^\dagger \quad \text{as} \quad \delta \downarrow 0.
\end{equation}
Some examples of regularizations methods are Tikhonov, Landweber, spectral cut-off, asymptotic regularization, approximate inverse and mollification (see, e.g. \cite{Alibaud,xavier_bonnefond,book_engl,book_kirsch,unified_app_reg_meth,paper_louis}).
As a matter of fact, we would like to get estimates on the error committed while approximating~$x^\dagger$ by~$x^\delta =  R_{\Lambda(\delta,y^\delta)} y^\delta$. 

It is well known that for arbitrary~$x^\dagger \in X$, the convergence of~$x^\delta$ towards~$x^\dagger$ is arbitrarily slow (see, e.g. \cite{book_engl,arbitrary_slow_conv}). But still, by allowing smoothness of the solution~$x^\dagger$, convergence rates could be established. Standard smoothness conditions known as H\"older type source condition take the form 
\begin{equation}
\label{holder source cond}
x^\dagger \in X_{\mu}(\rho) = \left\lbrace (T^*T)^\mu w, \quad w \in X \quad \text{s.t.}  \quad ||w|| \leq \rho \right\rbrace,
\end{equation}
where~$\mu$ and~$\rho$ are two positive constants. However such source conditions have shown their limitations as they are too restrictive in many problems and do not yield a natural interpretation. For this reason, general source conditions have been introduced in the following form:
\begin{equation}
\label{gen_sour_cond}
x^\dagger \in X_{\varphi}(\rho) = \left\lbrace \varphi(T^*T) w, \quad w \in X \quad \text{s.t.}  \quad ||w|| \leq \rho \right\rbrace,
\end{equation}
where~$\rho$ is a positive constant and $\varphi: [0, ||T^*T||] \to \Rb_+$ is an index function, i.e. a non-negative monotonically increasing continuous function  satisfying~$\varphi(\lambda) \to 0$ as~$\lambda \downarrow 0$. An interesting discussion on these source conditions can be found in \cite{how_gen_are_gen_sour_cond} where the author explores how general source conditions of the form \eqref{gen_sour_cond} are. Once the solution~$x^\dagger$ satisfies a smoothness condition i.e.~$x^\dagger$ belongs to a proper subspace~$M$ of~$X$, it is possible to derive convergence rates and the next challenge is about optimality. More precisely, for a regularization method $R: Y \to X$, we are interested in the worst case error:
\begin{equation}
\label{worst_case_err}
\Delta(\delta,R,M) := \sup \left\lbrace || R y^\delta - x^\dagger ||,\,\, x^\dagger \in M,\,\, y^\delta \in Y, \,\, \text{s.t.} \,\, ||y^\delta - T x^\dagger || \leq \delta \right\rbrace,
\end{equation}
and we would like a regularization which minimizes this worst case error. In this respect, a regularization method~$\bar{R}: Y \to X$ is said to be optimal if it achieves the minimum worst case error over all regularization methods, i.e. if 
\begin{equation*}
\Delta(\delta,\bar{R},M)  = \Delta(\delta,M) := \inf_R \Delta(\delta,R,M).
\end{equation*}
Similarly, a regularization is said to be order optimal if it achieves the minimum worst case error up to a constant greater than one, i.e. if 
\begin{equation*}
\Delta(\delta,\bar{R},M)  \leq C \Delta(\delta,M).
\end{equation*}
for some constant $C >1$.
When the subset $M$ is convex and balanced, it is shown in \cite{survey_on_optimal_recovery} that
\begin{equation}
\label{lower_bound_min_worst_case_error}
\omega(\delta,M)  \leq \Delta(\delta,M) \leq 2 \omega(\delta,M),
\end{equation}
where $\omega(\delta,M)$ is the modulus of continuity of the operator~$T$ over~$M$ i.e.
\begin{equation}
\label{modulus_continuity_T_over_M}
\omega(\delta,M) = \sup \left\lbrace ||x||, \,\, x \in M, \,\, \text{s.t.} \,\, || Tx|| \leq \delta \right\rbrace.
\end{equation}
In other words, we get the following: 
\begin{equation}
\label{min_worst_case_err_via_mud_cont}
\Delta(\delta,X_\varphi(\rho)) = \mathcal{O}\left( \omega(\delta,X_\varphi(\rho)) \right).
\end{equation}
Recall that, under mild assumptions on the index function~$\varphi$, the supremum defining the modulus of continuity is achieved and a simple expression of $\omega(\delta,X_\varphi(\rho))$ in term of function $\varphi$ is available (see, e.g. \cite{reg_exp_ipp,geom_lin_ipp,optimality_under_gen_sour_cond}). Let us remind that a relevant notion in the study of optimality of a regularization method is qualification. In fact, the qualification of a regularization measures the capability of the method to take into account smoothness assumptions on the solution~$x^\dagger$, i.e. the higher the qualification, the more the method is able to provide best rates for very smooth solutions.

Besides optimality, converse results and saturation results are also important aspects of regularization algorithms (see, \cite{book_engl,satur_reg_meth,on_conv_and_satur_results1,on_conv_and_satur_results2}). For converse results, we are interested in the following: given a particular convergence rate of $||x^\delta - x^\dagger||$ towards~$0$, which smoothness condition does the solution~$x^\dagger$ needs to satisfy? Saturation results are about the maximal smoothness on the solution~$x^\dagger$ for which a regularization method can still deliver the best rates of convergence. Finally, another significant aspect of regularization is the selection of the regularization parameter i.e. finding a function $\Lambda(\delta,y^\delta)$ which guarantees convergence and possibly order-optimality.

Coming back to \eqref{gen_sour_cond}, notice that a very interesting subclass of general source conditions are logarithmic source conditions expressed as:
\begin{equation}
\label{log_sour_cond}
x^\dagger \in X_{f_p}(\rho) = \left\lbrace (-\ln(T^*T))^{-p} w, \quad w \in X \quad \text{s.t.}  \quad ||w|| \leq \rho \right\rbrace,
 \end{equation}
where~$p$ and~$\rho$ are positive constants and $T$ satisfies $\norm{T^*T}_{} <1$. Such smoothness conditions have clear interpretations in term of Sobolev spaces in exponentially ill-posed problems (see, e.g. \cite{reg_exp_ipp,optimality_under_gen_sour_cond}). The latter class includes several problems of great importance such as backward heat equation, sideways heat equation, inverse problem in satellite gradiometry, control problem in heat equation, inverse scattering problems and many others (see, \cite{reg_exp_ipp}). Because of the importance of exponentially ill-posed problems, it is desirable to design regularization methods particularly suitable for this class of problems. It is precisely the aim of this paper to provide such a regularization scheme. 

 In the next section, we define the new regularization method using both the variational formulation and the definition in terms of the so called \textit{generator} function~$g_\alpha$. A brief comparison with the Tikhonov method is done. Moreover basic estimates on the \textit{generator} function $g_\alpha$ and its corresponding \textit{residual} function $r_\alpha$ are also carried out.
 
 Section \ref{section op_under_log_sour_cond} is devoted to optimality of the new method. Here we recall well known optimality results under general source conditions of the form \eqref{gen_sour_cond} (see, \cite{analysis_prof_funct,reg_exp_ipp,geom_lin_ipp,morozov_under_gen_sour_cond,
optimality_under_gen_sour_cond}). For the specific case of logarithmic source conditions, qualification of the method is given and order optimality is shown. % and converse results are also presented. 
Next we study optimality under general source conditions. %Moreover a posteriori parameter choice rule proposed in \cite{geom_lin_ipp} is presented.

 In Section \ref{section_comp_reg_meth}, we present a comparative analysis of the new method with Tikhonov method, spectral cut-off, asymptotic regularization and conjugate gradient.
 
 Section \ref{section_num_illust} is about numerical illustrations. In this section, in order to confirm our prediction of better performance of the new method compared to Tikhonov and spectral cut-off in instance of exponentially ill-posed problems, we numerically compare the efficiency of the five regularization methods on three test problems coming from literature: A problem of image reconstruction taken from \cite{shaw}, a Fredholm integral equation of the first kind found in \cite{baart} and an inverse heat equation problem. 
 
 Finally in Section \ref{Section_par_sel_rules}, for a fully applicability of the new method, we exhibit heuristic selection rules which fit with the new regularization technique. Moreover, we also compare the five regularization methods for each heuristic parameter choice rule under consideration.

\section{The new regularization method}{\label{section new_reg_meth}}
For the sake of simplicity, we assume henceforth that the operator $T$ is injective. Hereafter, we set a positive number $a$ such that the operator norm of $T^*T$ is less than $a$ i.e. $||T^*T||_{} \leq a$. In the sequel, we assume that $a<1$ which is always possible by scaling Equation \eqref{main equation}.

Let us consider the general variational formulation of a regularization method
\begin{equation}
\label{gen var formul}
x_\alpha = \argmin_{x \in X} \F(Tx,y) + \P(x,\alpha)
\end{equation}
where $\F(Tx,y)$ is the fit term, $\P(x,\alpha)$ is the penalty term and $\alpha>0$ is the regularization parameter. We recall that the fit term aims at fitting the model, the penalty term aims at introducing stability in the initial model $T x =y$ and the regularization parameter $\alpha$ controls the level of regularization.

In most cases, the fit term $\F(Tx,y)$ is nothing but
\begin{equation}
\label{fit term}
\F(Tx,y) = || Tx-y||_{}^2
\end{equation}
and the penalty term depends on the regularization method. For instance, for Tikhonov regularization, $\P(x,\alpha)$ is given by 
\begin{equation}
\label{pen term tik}
\P(x,\alpha) = \alpha ||x||_{}^2.
\end{equation}
This penalization can sometimes compromises the quality of the resulting approximate solution $x_\alpha$. Indeed, let $X = L
^2(\Rb^n)$, then by Parseval identity, we see that
\begin{equation}
\label{pen tik fourier}
\P(x,\alpha) = \alpha || \hat{x} ||_{L^2(\Rb^n)}^2
\end{equation}
where $\hat{x}$ is the Fourier transform of $x$. Equation \eqref{pen tik fourier} implies that the stability is introduced by uniformly penalizing all frequency components irrespective of the magnitude of frequencies. Yet, it is well known that the instability of the initial problem comes from high frequency components  on the contrary to low frequency components.% which do not generally cause any issue.

Let us introduce the following penalty term where the regularization parameter $\alpha$ is no more defined as a weight but as an exponent:
\begin{equation}
\label{def penalty nrm}
\P(x,\alpha) = \norm{\left[I - (T^*T)^{\sqrt{\alpha}} \right] x}_{}^2.
\end{equation}
In \eqref{def penalty nrm}, $(T^*T)^{\sqrt{\alpha}}$ is defined via the spectral family $(E_\lambda)_\lambda$ associated to the self-adjoint operator $T^*T$ i.e.
$$
(T^*T)^{\sqrt{\alpha}} x  = \int_{\lambda = 0}^{||T^*T||_+} \lambda^{\sqrt{\alpha}} \ud E_\lambda x.
$$
We keep the fit term defined in \eqref{fit term} and then the variational formulation of our new regularization method is given by 
\begin{equation}
\label{var for nrm}
x_\alpha = \argmin_{x \in X} || Tx -y||_{}^2 +  \norm{\left[I - (T^*T)^{\sqrt{\alpha}} \right] x}_{}^2.
\end{equation}
From the first order optimality condition, we get that $x_\alpha$ is the solution to the linear equation :
$$
\left[ T^*T +  \left(I - (T^*T)^{\sqrt{\alpha}} \right)^2 \right] x = T^* y,
$$
that is,
\begin{equation}
\label{formula app nrm}
x_\alpha = \left[ T^*T +  \left(I - (T^*T)^{\sqrt{\alpha}} \right)^2 \right]^{-1} T^* y.
\end{equation}
From \eqref{formula app nrm}, we see that the new method can also be defined via the so called \textit{generator} function $g_\alpha$, i.e.
\begin{equation}
\label{formula app nrm with gen func}
x_\alpha = g_\alpha(T^*T) T^* y,
\end{equation}
with the function $g_\alpha$ defined by 
\begin{equation}
\label{def_gen_func_new_method}
g_\alpha(\lambda) = \frac{1}{\lambda + (1-\lambda^{\sqrt{\alpha}})^2}, \quad \lambda \in (0, ||T^*T||].
\end{equation}
Let us also define the \textit{residual} function $r_\alpha$ corresponding to $g_\alpha$ as follows
\begin{equation}
\label{def_res_func_new_method}
r_\alpha(\lambda) := 1 - \lambda g_\alpha(\lambda) = \frac{(1-\lambda^{\sqrt{\alpha}})^2}{\lambda + (1-\lambda^{\sqrt{\alpha}})^2} , \quad \lambda \in (0, ||T^*T||].
\end{equation}
The functions $g_\alpha$ and $r_\alpha$ defined in \eqref{def_gen_func_new_method} and \eqref{def_res_func_new_method} are important since they will be repeatedly used in the convergence analysis of the regularization method. In fact, the regularization error $x^\dagger - x_\alpha$ and the propagated error $x_\alpha - x_\alpha^\delta$ are expressed via the functions $r_\alpha$ and $g_\alpha$ as follows:
\begin{equation*}
x^\dagger - x_\alpha = r_\alpha(T^*T)x^\dagger,\quad x_\alpha - x_\alpha^\delta = g_\alpha(T^*T)T^*(y - y^\delta).
\end{equation*}
Finally, notice that the function $g_\alpha$ defined in \eqref{def_gen_func_new_method} indeed satisfies the basic requirements for defining a regularization method i.e. 
\begin{itemize}
\item[a)] $g_\alpha$ is continuous,
\item[b)] $ \forall\, \alpha>0, \quad \sup_{\lambda \in (0,||T^*T||]} \lambda g_\alpha(\lambda) \leq 1 < \infty$,
\item[c)] $\lim_{\alpha \downarrow 0} \,g_\alpha(\lambda) = 1/\lambda.$
\end{itemize}

From b) and c), we deduce the convergence of the new regularization method by application of \cite[Theorem 4.1]{book_engl}. Before going to optimality results, let us state some basic estimates (proven in the appendix) about the functions $g_\alpha$ and $r_\alpha$.
\begin{Proposition}
\label{prop 2} 
Let the function $g_\alpha$ be defined by \eqref{def_gen_func_new_method}. Then for all $a <1$ and $\alpha <1$, 
\begin{equation}
\label{bound fct prop err}
\sup_{\lambda \in (0,a]} \sqrt{\lambda} g_\alpha(\lambda) = \bigO{\frac{1}{\sqrt{\alpha}}}.
\end{equation}
\end{Proposition}

\begin{Lemma}
\label{Lemma bound resid_func}
For all $\alpha$ and $\lambda$ satisfying $0< \alpha \leq \lambda  < 1$, the following estimates hold for the function $r_\alpha$ defined in \eqref{def_res_func_new_method}:
\begin{equation}
\label{bound_res_func_new_reg_meth}
r_\alpha(\lambda) \leq \frac{9}{4} \left( \frac{\alpha |\ln(\lambda)|^2}{\lambda +\alpha |\ln(\lambda)|^2}\right).
\end{equation}
\end{Lemma}

\section{Optimality results}{\label{section op_under_log_sour_cond}}
Before studying the optimality of the method presented in Section \ref{section new_reg_meth}, we need first to recall general optimality results under source condition of the form \eqref{gen_sour_cond}. For doing so, let us specify assumptions on the function $\varphi$ which defines the source set $X_\varphi(\rho)$. 
\begin{Assumption}
\label{Assumption_index_function}
The function $\varphi: (0,a] \to \Rb_+$ is continuous, monotonically increasing and satisfies:
\begin{itemize}
\item[(i)] $\lim_{\lambda \downarrow 0} \varphi(\lambda) = 0$,
\item[(ii)] the function $\phi: (0,\varphi^2(a)] \to (0,a\varphi^2(a)]$ defined by 
\begin{equation}
\label{def_func_phi}
\phi(\lambda) = \lambda (\varphi^2)^{-1}(\lambda)
\end{equation}
is convex.
\end{itemize}
\end{Assumption}

Under Asummption \ref{Assumption_index_function} on the function $\varphi$, the following result from \cite{optimality_under_gen_sour_cond} holds and we can then defines optimality under source condition \eqref{gen_sour_cond}.
\begin{Theorem}
\label{theorem optimality}
Let $X_\varphi(\rho)$ be as in \eqref{gen_sour_cond} and let Assumption \ref{Assumption_index_function} be fulfilled. Let the function $\phi$ be defined by \eqref{def_func_phi}. Then
\begin{equation}
\label{bound_mudulos_of_cont}
\omega(\delta,X_\varphi(\rho)) \leq \rho \sqrt{\phi^{-1}\left(\frac{\delta^2}{\rho^2} \right)}.
\end{equation}
Moreover, if $\delta^2/\rho^2 \in \sigma\left(T^*T \varphi^2(T^*T) \right)$, then equality holds in \eqref{bound_mudulos_of_cont}.
\end{Theorem}
A similar result to this theorem can be found in \cite[Section 2]{reg_exp_ipp}, and \cite[Section 3]{geom_lin_ipp}. 
\begin{Remark}
\label{remark optim}
In \cite{geom_lin_ipp}, the results corresponding to Theorem \ref{theorem optimality} are given in term of the function $\Theta: (0,a] \to (0,a\varphi(a)]$ defined by:
\begin{equation}
\label{def_func_Theta}
\Theta(\lambda) = \sqrt{\lambda}\varphi(\lambda).
\end{equation}
Then, by simple computations, we can find that
\begin{equation}
\label{equality_best_bound_from_two_ref}
\rho \,\sqrt{\phi^{-1}\left(\frac{\delta^2}{\rho^2}\right)} = \rho \,\varphi \left(\Theta^{-1}(\delta/\rho) \right).
\end{equation}
In such a case, the convexity of the function~$\phi$ defined in \eqref{def_func_phi} is equivalent to the convexity of the function $\chi(\lambda) = \Theta^2\left( (\varphi^2)^{-1}(\lambda)\right)$ and the condition $\delta^2/\rho^2 \in \sigma\left(T^*T \varphi^2(T^*T) \right)$ which allows to get the equality in \eqref{bound_mudulos_of_cont} is equivalent to $\delta/\rho \in \sigma \left( \Theta(T^*T)\right)$.
\end{Remark}
From Theorem \ref{theorem optimality} and Remark \ref{remark optim}, we can deduce that under the source condition \eqref{gen_sour_cond} and Assumption~\ref{Assumption_index_function}, the best possible worst case error is $\rho \,\varphi \left(\Theta^{-1}(\delta/\rho) \right)$ whence the following definition.
\begin{Definition}[Optimality under general source conditions]
\label{Def_opt_under_gen_sour_cond}
Let Assumption \ref{Assumption_index_function} be satisfied and consider the source condition $x^\dagger \in X_\varphi(\rho)$. A regularization method $R(\delta): Y \to X$ is said to be:
\begin{itemize}
\item optimal if $\Delta(\delta,R(\delta),X_{\varphi}(\rho)) \leq \rho \,\varphi \left(\Theta^{-1}(\delta/\rho) \right)$;
\item order optimal if $\Delta(\delta,R(\delta),X_{\varphi}(\rho)) \leq C \rho \,\varphi \left(\Theta^{-1}(\delta/\rho) \right)$
for some constant $C \geq 1$;
\item quasi-order optimal if for all $\epsilon >0$, $\Delta(\delta,R(\delta),X_{\varphi}(\rho)) = \mathcal{O}\left( f_\epsilon(\delta) \right) $ where the function $f_\epsilon: \Rb_+ \to \Rb_+$ converges to $\varphi \left(\Theta^{-1}(\delta/\rho) \right)$ as~$\epsilon$ decreases to $0$ i.e. for all $\delta >0$, $f_\epsilon(\delta) \to \varphi \left(\Theta^{-1}(\delta/\rho) \right)$ as~$\epsilon$ decreases to $0$.

\end{itemize}
\end{Definition}

Having defined the optimality under general source conditions, let us now consider the particular case of logarithmic source conditions. For logarithmic source conditions, the function $\varphi$ equals the function $f_p: (0,a] \to \Rb_+$ defined by:
\begin{equation}
\label{def_func_f_p}
f_p(\lambda) = (-\ln(\lambda))^{-p}.
\end{equation}
Next it is easy to see that the only point to check in Assumption \ref{Assumption_index_function} is the convexity of the function $\phi$ defined in \eqref{def_func_phi}. Precisely, for the index function $f_p$, this function is $\phi_p: (0, \ln(1/a)^{-2p}] \to (0, a \ln(1/a)^{-2p}]$ defined by 
\begin{equation*}
%\label{function_phi_p}
\phi_p(\lambda) = \lambda\exp(-\lambda^{-1/2p})
\end{equation*}
 which was proven to be convex on the interval $[0,1]$ in \cite{tik_for_fini_infini_smooth_op}. In order to fulfill Assumption \ref{Assumption_index_function} and avoid the singularity of the function $f_p$ at $\lambda = 1$, we assume that $a \leq \exp(-1) < 1$, i.e.~$||T^*T||_{} \leq \exp(-1)$. Notice that this is not actually a restriction, since Equation \eqref{main equation} can always be rescaled in order to meet this criterion. 

Due to \eqref{bound_mudulos_of_cont} it suffices to compute $\sqrt{\phi_p^{-1}\left(\delta^2/\rho^2\right)}$ in order to define the optimality in logarithmic source conditions. Thanks again to \cite{tik_for_fini_infini_smooth_op}, we have that
\begin{equation}
\label{def_inv_func_phi_p}
\sqrt{\phi_p^{-1}(s)} = f_p(s) (1 + o(1)) \,\, \text{as} \,\, s \to 0.
\end{equation}

Hence, we deduce the following definition of optimality in case of logarithmic source condition.

\begin{Definition}[Optimality under logarithmic source condition]
\label{Def opt log sour cond}
Consider logarithmic source condition \eqref{log_sour_cond}, on defining $f_p$ as in \eqref{def_func_f_p}, a regularization method $R(\delta): Y \to X$ is said to be:
\begin{itemize}
\item optimal if $\Delta(\delta,R(\delta),X_{f_p}(\rho)) \leq \rho f_p(\delta^2/\rho^2) (1 + o(1)) \,\, \text{as} \,\, \delta \to 0$,
\item order optimal if $\Delta(\delta,R(\delta),X_{f_p}(\rho)) \leq C \rho f_p(\delta^2/\rho^2) (1 + o(1)) \,\, \text{as} \,\, \delta \to 0 \,$.
\end{itemize}
%where the function $f_p$ is defined in \eqref{def_func_f_p}.
\end{Definition}
In the sequel, we are interested to optimality with respect to the noise level $\delta$. In this respect, we can characterize the order-optimality under logarithmic source conditions as follows.
\begin{Remark}
By definition of the function $f_p$, we get that $\bigO{f_p(\delta^2/\rho^2)} = \bigO{f_p(\delta)}$ as $\delta \to 0$. Hence, equivalently to Definition \ref{Def opt log sour cond}, a regularization method $R(\delta): Y \to X$ is said to be order optimal under logarithmic source condition if 
$$
\Delta(\delta,R(\delta),X_{f_p}(\rho)) = \bigO{f_p(\delta)} \quad \textrm{as} \quad \delta \to 0.
$$
%Similarly, a regularization method $R(\delta): Y \to X$ is said to be order optimal under general source condition $x^\dagger \in X_\varphi(\rho)$ if 
%$$
%\Delta(\delta,R(\delta),X_{f_p}(\rho)) = \bigO{\varphi(\Theta^{-1}(\delta)} \quad \textrm{as} \quad \delta \to 0.
%$$
\end{Remark}

\subsection{Optimality under logarithmic source conditions}
Having given all the necessary definitions, let us now study the optimality of the  method proposed in Section \ref{section new_reg_meth}. %Let us state the following qualification result on the regularization method.
\begin{Proposition}
\label{Proposition qualification}
The regularization $g_\alpha$ defined by \eqref{def_gen_func_new_method} has qualification $f_p$. That is:
\begin{equation}
\label{qualif_new_meth}
\sup_{0<\lambda \leq a} \, |r_\alpha(\lambda)|\,f_p(\lambda)  = \mathcal{O}\left( f_p(\alpha) \right).
\end{equation}
\end{Proposition}
The proof the Proposition \ref{Proposition qualification} heavily relies on the following lemma which is proven in the appendix.
\begin{Lemma}
\label{Lemma 2}
Let $p$ and $\alpha$ be two positive numbers with $\alpha \leq \bar{\alpha} <1$, let $a \in (0,1)$ and $\Psi_{p,\alpha}: (0,a] \to \Rb_+$ be the function defined by
\begin{equation}
\label{def_func_chi}
\Psi_{p,\alpha}(\lambda) = \frac{|\ln(\lambda)|^{2-p}}{\lambda + \alpha |\ln(\lambda)|^2}.
\end{equation}
Then, the following hold:
\begin{itemize}
\item[(i)] The function $\Psi_{p,\alpha}$ is well defined and differentiable on $(0,a]$, and its derivative is given by
\begin{equation}
\label{der_func_chi}
\Psi_{p,\alpha}'(\lambda) = \frac{\lambda^{-1}|\ln(\lambda)|^{1-p}}{(\lambda+\alpha |\ln(\lambda)|^2 )^2}  h(\lambda),
\end{equation}
where
\begin{equation}
\label{def_func_h_der_func_chi}
h(\lambda) =  \alpha p |\ln(\lambda)|^2 - \lambda \left(2-p + |\ln(\lambda)|\right).
\end{equation}
\item[(ii)] If $p \leq 2$, there exists at least one $\lambda(\alpha,p)$ where $h$ vanishes. Moreover for every such $\lambda(\alpha,p)$, the following holds
\begin{equation}
\label{estimate root of h}
\lambda(\alpha,p) \simeq \alpha |\ln(\alpha)|,
\end{equation}
that is, there exists two constants $c_1$ and $c_2$ depending on $p$ only such that
$$
 c_1 \alpha |\ln(\alpha)| \leq \lambda(\alpha,p) \leq c_2 \alpha |\ln(\alpha)|.
$$
Moreover, this result still holds if $p>2$, $\lambda < c \leq \exp(2-p)$ and $\alpha$ is small.
\item[(iii)] The supremum of the function $\Psi_{p,\alpha}$ on $(0,a]$ satisfies
\begin{equation}
\label{supremum_func_chi}
\sup_{0 < \lambda \leq a} \Psi_{p,\alpha}(\lambda) = \mathcal{O}\left( \alpha^{-1} |\ln(\alpha)|^{-p} \right).
\end{equation}
\end{itemize}
\end{Lemma}
Having stated the above lemma, the proof of Proposition \ref{Proposition qualification} easily follows:

\begin{Proof}
If $\lambda \leq \alpha$ then the monotonicity of the function $f_p$ and the fact that the residual function $r_\alpha$ is bounded by 1 on $(0,a]$ yields \eqref{qualif_new_meth}. If $\lambda \geq \alpha$
then from Lemma \ref{Lemma 2}, we deduce that 
$$
\sup_{0< \lambda \leq a} \frac{\alpha |\ln(\lambda)|^2}{\lambda +\alpha |\ln(\lambda)|^2} f_p(\lambda) = \mathcal{O}\left( f_p(\alpha) \right)
$$
which together with Lemma \ref{Lemma bound resid_func} yields \eqref{qualif_new_meth}.
\end{Proof}

From Proposition \ref{Proposition qualification}, we deduce the following optimality result.
\begin{Theorem}
\label{Theorem opt under log sour cond}
Let $p>0$, $x^\dagger \in X_{f_p}(\rho)$, and $y^\delta \in Y$ satisfying \eqref{noise level equation} with $y = T x^\dagger$. Assume that $||T^*T||_{} \leq \exp(-1)$ and let $x(\delta) = g_{\alpha(\delta)}(T^*T)T^*y^\delta$ with the function $g_\alpha$ being defined by \eqref{def_gen_func_new_method} and let $\alpha(\delta) = \Theta_p^{-1}(\delta)$
with $\Theta_p$ defined by
\begin{equation}
\label{def_func_Theta_p}
\Theta_p(\lambda) = \sqrt{\lambda}(\ln(1/\lambda))^{-p}.
\end{equation}
Then the order optimal estimate
\begin{equation}
\label{opt rate log sour cond}
|| x^\dagger - x(\delta) || = \mathcal{O} \left( f_p(\delta) \right) \quad \text{as} \quad \delta \to 0 
\end{equation}
holds.
Thus the regularization $g_\alpha$ defined by \eqref{def_gen_func_new_method} is order optimal under logarithmic source conditions.
\end{Theorem}

\begin{Proof}%[of Theorem \ref{Theorem opt under log sour cond}]
As usual, we start with the following splitting
\begin{equation}
\label{splitting error }
||x^\dagger - x_\alpha^\delta || \leq ||x^\dagger - x_\alpha||  + || x_\alpha - x_\alpha^\delta||.
\end{equation}
Using that $x^\dagger - x_\alpha = r_\alpha(T^*T) x^\dagger$, $x_\alpha - x_\alpha^\delta = g_\alpha(T^*T)T^*(y -y^\delta)$ together with the source condition $x^\dagger \in X_{f_p(\rho)}$, we deduce that:
\begin{equation}
\label{bound reg error}
||x^\dagger - x_\alpha|| \leq C_1 \sup_{\lambda \in (0,a]} r_\alpha(\lambda)\, f_p(\lambda)
\end{equation}
and
\begin{equation}
\label{bound prop error}
|| x_\alpha - x_\alpha^\delta|| \leq \delta \,\, C_2 \sup_{\lambda \in (0,a]} \sqrt{\lambda} g_\alpha(\lambda).
\end{equation}

By applying the propositions \ref{prop 2} and \ref{Proposition qualification} to \eqref{bound reg error}, \eqref{bound prop error} and using \eqref{splitting error }, we get that 
\begin{equation}
\label{final bound err log sour cond}
||x^\dagger - x_\alpha^\delta || \leq C_1' f_p(\alpha) + C_2' \frac{\delta}{\sqrt{\alpha}},
\end{equation}
where $C_1'$ and $C_2'$ are constants independent of $\alpha$ and $\lambda$.
Hence, by taking $\alpha:= \Theta_p^{-1}(\delta)$, the estimate in \eqref{opt rate log sour cond} follows from
\begin{equation*}
||x^\dagger - x(\delta)|| = \bigO{f_p(\Theta_p^{-1}(\delta))} = \bigO{f_p(\delta^2)} = \bigO{f_p(\delta)}.
\end{equation*}
\end{Proof}

\begin{Corollary}
Let $p>0$, $x^\dagger \in X_{f_p}(\rho)$, and $y^\delta \in Y$ satisfying \eqref{noise level equation} with $y = T x^\dagger$. Assume that $||T^*T||_{} \leq \exp(-1)$ and let $x(\delta) = g_{\alpha(\delta)}(T^*T)T^*y^\delta$ with the function $g_\alpha$ being defined by \eqref{def_gen_func_new_method} and  $\alpha(\delta) = \delta$. 
Then the order optimal estimate 
$$
|| x^\dagger - x(\delta) || = \mathcal{O} \left( f_p(\delta) \right) \quad \text{as} \quad  \delta \to 0
$$
holds. Thus the regularization $g_\alpha$ defined by \eqref{def_gen_func_new_method} is order optimal under logarithmic source conditions with an a-priori parameter choice rule independent of the smoothness of the solution $x^\dagger$.
\end{Corollary}

\begin{Proof}
By considering $\alpha(\delta) = \delta$ in \eqref{final bound err log sour cond}, we get 
\begin{equation*}
||x^\dagger - x_\alpha^\delta || \leq C_1' f_p(\delta) + C_2' \sqrt{\delta} = \bigO{f_p(\delta)} \quad \textrm{as} \quad \delta \to 0,
\end{equation*}
since $\sqrt{\delta} = \bigO{f_p(\delta)}$ as $\delta \to 0$.
\end{Proof}

The next proposition describes a Morozov-like discrepancy rule which leads to order-optimal convergence rates under logarithmic source conditions.
\begin{Proposition}
\label{Prop order optimal conv rates under posteriori rule}
Let $p>0$, $x^\dagger \in X_{f_p}(\rho)$, and $y^\delta \in Y$ satisfying \eqref{noise level equation} with $y = T x^\dagger$. Assume that $||T^*T||_{} \leq \exp(-1)$ and consider the a-posteriori parameter choice rule
\begin{equation}
\label{morozov_like_par_sel_rule}
\alpha(\delta,y^\delta) = \sup \left\lbrace \alpha>0, \quad ||T x_\alpha^\delta - y^\delta|| \leq \delta + \sqrt{\delta}\right\rbrace.
\end{equation}
Let $x(\delta) = g_{\alpha(\delta,y^\delta)}(T^*T)T^*y^\delta$ with the function $g_\alpha$ defined by \eqref{def_gen_func_new_method}, then the order optimal estimate
\begin{equation}
\label{eq tnnggggoifng}
|| x^\dagger - x(\delta) || = \mathcal{O} \left( f_p(\delta) \right) \quad \text{as} \quad \delta \to 0
\end{equation}
holds.
Thus the regularization $g_\alpha$ defined by \eqref{def_gen_func_new_method} is order optimal under logarithmic source conditions with the a-posteriori parameter choice rule defined by \eqref{morozov_like_par_sel_rule}.
\end{Proposition}
The proof of Proposition \ref{Prop order optimal conv rates under posteriori rule} is deferred to Appendix. 
%From \eqref{final bound err log sour cond}, we can define a simple a-posteriori parameter selection rule $\alpha(\delta,y^\delta)$ which also leads to order optimal convergence rates under logarithmic source conditions.
%\begin{Proposition}
%\label{Prop order optimal conv rate simple post rule}
%Consider the setting of Proposition \ref{Prop order optimal conv rates under posteriori rule} except that we consider
%\begin{equation}
%\label{morozov_simple_par_sel_rule}
%\alpha(\delta,y^\delta) = \sup \left\lbrace \alpha \in (0,\delta), \quad ||T x_\alpha^\delta - y^\delta|| \leq \delta \right\rbrace.
%\end{equation}
%Let $x(\delta) = g_{\alpha(\delta,y^\delta)}(T^*T)T^*y^\delta$ with the function $g_\alpha$ defined by \eqref{def_gen_func_new_method}, then the order optimal estimate
%\begin{equation*}
%|| x^\dagger - x(\delta) || = \mathcal{O} \left( f_p(\delta) \right) \quad \text{as} \quad \delta \to 0
%\end{equation*}
%holds.
%\end{Proposition}
%\begin{Proof}
%From \eqref{morozov_simple_par_sel_rule}, we see that $\alpha(\delta,y^\delta) = \tau \delta$ with $\tau \leq 1$. Hence from \eqref{final bound err log sour cond}, we get that 
%$$
%|| x^\dagger - x(\delta)|| \leq C_1' f_p(\tau \delta)  + \frac{C_2'}{\sqrt{\tau}} \sqrt{\delta} = \bigO{f_p(\delta)} \quad \text{as} \quad \delta \to 0,
%$$
%since $f_p(\tau \delta) = \bigO{f_p(\delta)}$ and $\sqrt{\delta} =  \bigO{f_p(\delta)}$ as $\delta \to 0$.
%\end{Proof}

\begin{flushleft}
\textbf{A converse result}
\end{flushleft}
Theorem \ref{Theorem opt under log sour cond} establishes that the logarithmic source condition \eqref{log_sour_cond} is sufficient to imply the rate $f_p(\delta)$ in \eqref{opt rate log sour cond}. Now we are going to prove that the logarithmic source condition \eqref{log_sour_cond} is not only sufficient but also almost necessary. % Let us start with the noise free case.
 The following result based on \cite[Theorem 8]{reg_exp_ipp} establishes a converse result in the noise free case for the new regularization method.
\begin{Theorem}
\label{Theo conv result noise free case}
Let $x_\alpha = g_\alpha(T^*T)Ty$ with $y = T x^\dagger$ and let the function $g_\alpha$ be defined in \eqref{def_gen_func_new_method}. Then the estimate
\begin{equation}
\label{rate noise free case log}
|| x^\dagger - x_\alpha|| = \bigO{f_p(\alpha)}
\end{equation}
implies that $x^\dagger \in X_{f_q}(\rho)$ for some $\rho >0$ for all $ 0 < q <p$.
\end{Theorem}
The proof consists in checking that the function $g_\alpha$ defined in \eqref{def_gen_func_new_method} satisfies all the conditions stated in Theorem 8 of \cite{reg_exp_ipp}. More precisely, we just need to check that there exists a constant $C_g>0$ such that
$$
\sup_{\lambda \in (0,||T^*T||]} g_\alpha(\lambda)  \leq \frac{C_g}{\alpha}.
$$
But, from \eqref{bound gen func nrm}, we see that the latter condition is obviously fulfilled.
%Now let us turn to the noisy case. The idea for the noisy case is to remark that the rate \eqref{rate noise free case log} corresponding to the noise free case is actually the same rate in \eqref{opt rate log sour cond} for the noisy case except that the $\alpha$ is replaced by $\delta$.
%
%\begin{Theorem}
%\label{Theo conv result noisy case}
%Let $x_\alpha^\delta = g_\alpha(T^*T)T y^\delta$ with the function $g_\alpha$ defined in \eqref{def_gen_func_new_method}. Then the estimates 
%\begin{equation}
%\label{rate noisy case log}
%\sup \left\lbrace \inf_{\alpha >0} || x^\dagger - x_\alpha^\delta||\quad  \Big\vert \quad || T x^\dagger - y^\delta|| \leq \delta \right\rbrace = \bigO{f_p(\delta)}
%\end{equation}
%implies that $x^\dagger \in X_{f_q}(\rho)$ for some $\rho >0$ for all $ 0 < q <p$.
%\end{Theorem}
%\begin{Proof}
%From \eqref{rate noisy case log}, by taking $y^\delta = y = T x^\dagger$, we get that 
%$$
%\inf_{\alpha >0} || x^\dagger - x_\alpha|| = \bigO{f_p(\delta)}
%$$
%where $x_\alpha = g_\alpha(T^*T)T^*y$. Next by letting $\hat{\alpha}$ be the minimizer of $|| x^\dagger - x_\alpha||$ and $\delta$ be equal to $ \hat{\alpha}$, we get
%$$
%|| x^\dagger - x_{\hat{\alpha}}|| =  \bigO{f_p(\hat{\alpha})}
%$$ 
%which from Theorem \ref{Theo conv result noise free case} implies that $x^\dagger \in X_{f_q}(\rho)$ for some $\rho >0$ for all $ 0 < q <p$.
%\end{Proof}

\subsection{Optimality under general source conditions}
Let us state the following quasi-optimal result under general source conditions.
\begin{Theorem}
\label{Theorem opt gen sour cond}
Let $p>0$, $x^\dagger \in X_{\varphi}(\rho)$, where $\varphi$ is a concave index function satisfying Assumption \ref{Assumption_index_function} and $y^\delta \in Y$ satisfying $||y - y^\delta|| \leq \delta$ with $y = T x^\dagger$ and $\delta \leq \Theta(a)$. Assume that $||T^*T||_{} \leq a\leq \exp(-1)$ and let $x(\delta) = g_{\alpha(\delta)}(T^*T)T^*y^\delta$ with the function $g_\alpha$ defined in \eqref{def_gen_func_new_method}.  For small positive $\epsilon$, let $\alpha(\delta) = \Theta_\epsilon^{-1}(\delta)$ where the function $\Theta_\epsilon$ is defined by
$
\Theta_\epsilon(\lambda) = \lambda^{-\epsilon} \Theta(\lambda)
$
with $\Theta$ given in \eqref{def_func_Theta}.

Then the estimate
\begin{equation*}
|| x^\dagger - x(\delta) || = \mathcal{O} \left( (\Theta_\epsilon^{-1}(\delta))^{-\epsilon} \varphi (\Theta_\epsilon^{-1}(\delta)) \right) \quad \text{as} \quad \delta \to 0
\end{equation*}
holds. 
Moreover, as $\epsilon \downarrow 0$, $(\Theta_\epsilon^{-1}(\delta))^{-\epsilon} \varphi (\Theta_\epsilon^{-1}(\delta)) \to \varphi (\Theta^{-1}(\delta))$.
Thus the regularization method defined via the function $g_\alpha$ given in \eqref{def_gen_func_new_method} is quasi-order optimal under general source conditions.
\end{Theorem}

\begin{Proof}
We study two cases: $\alpha \geq \lambda$ and $\alpha < \lambda$.
In the first case, $\sup_{(0,\exp{(-1)}]} r_\alpha(\lambda) \varphi(\lambda) \leq  \varphi(\alpha)$ by monotonicity of the function $\varphi$ and the order-optimality follows trivially. Let us study the main case when $\alpha < \lambda$.
From Lemma \ref{Lemma bound resid_func}, we get, for $\lambda \in (0,a]$,
\begin{eqnarray}
\label{eq bound.}
r_\alpha(\lambda) \varphi(\lambda) 
& \leq & \frac{9}{4} |\ln(\lambda)| ^2 \frac{\alpha}{\lambda + \alpha |\ln(\lambda)| ^2} \varphi(\lambda) \nonumber \\
& \leq & \frac{9}{4} |\ln(\alpha)| ^2 \frac{\alpha}{\lambda + \alpha |\ln(a)| ^2} \varphi(\lambda)  \nonumber \\
& \leq & \frac{9}{4} \alpha^{-\epsilon} (\alpha^{\epsilon/2} |\ln(\alpha)| )^2 \frac{\alpha}{\lambda + \alpha |\ln(a)| ^2} \varphi(\lambda)   \nonumber \\
& \leq & \frac{9}{4} \frac{4}{\epsilon^2} \alpha^{-\epsilon} \frac{\alpha \lambda}{\lambda + \alpha |\ln(a)| ^2} \frac{\varphi(\lambda)}{\lambda}  \nonumber \\
& \leq & \frac{9}{4} \frac{4}{\epsilon^2} \alpha^{-\epsilon} \frac{\alpha \lambda}{\lambda + \alpha |\ln(a)| ^2} \frac{\varphi(\alpha)}{\alpha} \quad \text{by concavity of} \,\, \varphi  \nonumber \\
& \leq & C_\epsilon  \alpha^{-\epsilon} \varphi(\alpha) .
\end{eqnarray}
Hence  
$
\sup_{(0,a]} r_\alpha(\lambda) \varphi(\lambda) \leq C_\epsilon \alpha^{-\epsilon} \varphi(\alpha).
$
From \eqref{bound reg error} and \eqref{bound prop error}, and \eqref{bound fct prop err} we get
$$
|| x^\dagger - x_\alpha^\delta|| \leq C_\epsilon \alpha^{-\epsilon} \varphi(\alpha) + \frac{\delta}{\sqrt{\alpha}}.
$$
By taking $\alpha(\delta) = \Theta_\epsilon^{-1}(\delta)$ with $\Theta_\epsilon(\lambda) = \lambda^{1/2 - \epsilon} \varphi(\lambda)$, we get
$$
|| x^\dagger - x(\delta) || = \mathcal{O} \left( (\Theta_\epsilon^{-1}(\delta))^{-\epsilon} \varphi (\Theta_\epsilon^{-1}(\delta)) \right).
$$
Now, it remains to show that $(\Theta_\epsilon^{-1}(\delta))^{-\epsilon} \varphi (\Theta_\epsilon^{-1}(\delta))$ converges to the optimal rate $\varphi(\Theta^{-1}(\delta))$ as $\epsilon$ goes to $0$. Let $\alpha_*  = \Theta^{-1}(\delta)$ and $\alpha_\epsilon = \Theta_\epsilon^{-1}(\delta)$, let us show that $\alpha_\epsilon$  converges to $\alpha_*$ as $\epsilon$ goes to $0$.
By the monotonicity of $\Theta_\epsilon$ for $\epsilon \in (0,1/2)$ and the fact that $\delta \leq \Theta(a)$ and $a<1$, we get that, for all $\epsilon \in (0,1/2)$,
$$
\frac{\delta}{\Theta(a)} \leq 1 < a^{-\epsilon}\quad  \Rightarrow \quad \delta \leq a^{-\epsilon}\Theta(a) = \Theta_{\epsilon}(a)\quad \Rightarrow \quad \alpha_\epsilon  =\Theta_{\epsilon}^{-1}(\delta) \leq a.
$$
Hence $\alpha_\epsilon \in (0,a]$ and the sequence $(\alpha_\epsilon)_\epsilon$ is bounded and thus it admits a converging subsequence. Let $(\alpha_{\epsilon_n})_n$ a converging subsequence of $(\alpha_\epsilon)_{\epsilon}$, and let $\tilde{\alpha}$ be its limit. Let us show that $\tilde{\alpha} = \alpha_*$.

Since $\alpha_{\epsilon_n} \to \tilde{\alpha}$ and $\Theta$ is continuous,
$\Theta(\alpha_{\epsilon_n}) \to \Theta(\tilde{\alpha})$. But 
$ \Theta(\alpha_{\epsilon_n}) = \alpha_{\epsilon_n}^{\epsilon_n} \Theta(\alpha*)$ since $\delta = \Theta(\alpha_*) $ and $\delta = \Theta_\epsilon(\alpha_\epsilon)$ for all small positive $\epsilon$.
So we get  
\begin{equation}
\label{ggggg}
\alpha_{\epsilon_n}^{\epsilon_n} \Theta(\alpha*) \to  \Theta(\tilde{\alpha}) \quad \rm{i.e.} \quad \alpha_{\epsilon_n}^{\epsilon_n} \to \frac{\Theta(\tilde{\alpha})}{\Theta(\alpha*)}.
\end{equation}
By the convergence of the sequence  $(\alpha_{\epsilon_n})_n$, we get that $\alpha_{\epsilon_n}^{\epsilon_n} = \exp{(\epsilon_n \ln (\alpha_{\epsilon_n}))}$ converges to $1$, \eqref{ggggg} proves that $\Theta(\tilde{\alpha})= \Theta(\alpha*)$ and by bijectivity of the function $\Theta$, we deduce that $\tilde{\alpha} = \alpha_*$.
Since the sequence $(\epsilon_n)_n$ was arbitrarily chosen, we deduce that the whole sequence $(\alpha_\epsilon)_\epsilon$ converges to $\alpha_*$ as $\epsilon \downarrow0$. 
Thus we deduce that $\alpha_\epsilon^{-\epsilon} \to 1$ and $\varphi(\alpha_\epsilon) \to \varphi(\alpha_*)$
which implies that
$$
(\Theta_\epsilon^{-1}(\delta))^{-\epsilon} \varphi (\Theta_\epsilon^{-1}(\delta)) \to \varphi(\Theta^{-1}(\delta)).
$$
\end{Proof}

For Holder type source conditions, Theorem \ref{Theorem opt gen sour cond} reduces to the following theorem.
\begin{Theorem}
\label{Theorem variant}
Consider the setting of Theorem \ref{Theorem opt gen sour cond} with the function $\varphi (t) = t^\mu$ i.e. $x^\dagger \in \ran{(T^*T)^\mu}$, then there exists an a priori selection rule $\alpha(\delta)$ such that the following holds:
\begin{equation}
\label{New convergence rate second scheme}
\norm{x^{\dagger} - x_{\alpha(\delta)}^{\delta}} = \begin{cases}
\mathcal{O}\left( \delta^{\frac{2 \sigma }{2 \sigma +1}}\right) \quad \forall\, \sigma < \mu , \,\, \text{if} \,\, \mu \leq 1  \\
\mathcal{O}\left( \delta^{\frac{2}{3}} \right) \quad , \,\, \text{if} \,\, \mu > 1.
\end{cases}
\end{equation}
\end{Theorem}

\begin{Remark}
By defining a variant of the new regularization method where the approximate solution~ $x_\alpha^\delta$ is defined as the solution of the optimization  problem
$$
x_\alpha^{\delta} = \argmin_{x \in X} || (T^*T)^{\sqrt{\alpha}} y^\delta - T x||_{}^2 + ||\left[I - (T^*T)^{\sqrt{\alpha}} \right] x ||_{}^2,
$$
we can prove order optimal rate under Holder type source condition but with a lower qualification index $\mu_0 = 1/2$. This variant is motivated by the mollification regularization method, where a target object defined as a smooth version of $x^\dagger$ is fixed prior to the regularization (see e.g. \cite{Alibaud,xavier_bonnefond}). In this respect, the target object here is given as $(T^*T)^{\sqrt{\alpha}} x^\dagger$. This choice is legitimated by the smoothness property of the operator $T$ and the fact that as $\alpha$ goes to $0$, this target object converges to the solution $x^\dagger$. The study of this variant and the corresponding optimality results is beyond the scope of this paper.
\end{Remark}

\section{A framework for comparison}{\label{section_comp_reg_meth}}

In the sequel, we are going to compare the new method with three continuous regularization methods: Tikhonov \cite{tikhonov}, spectral cut-off \cite{book_engl}, Showalter \cite{book_engl} and one iterative regularization method: conjugate gradient \cite{book_engl,book_kirsch}. We recall that the first three methods (Tikhonov, spectral cut-off and Showalter) are linear methods on the contrary to conjugate gradient which is an iterative non-linear regularization method. Obviously the new method, Tikhonov, spectral cut-off and Showalter are members of the family of general regularization methods defined via a \text{\it generator} function. Roughly speaking, each regularization method is defined via a so-called \textit{generator} function $g_\alpha^{reg}(\lambda)$ which converges pointwise to $1/\lambda$ as $\alpha$ goes to $0$ and the regularized solution $x_{\alpha,reg}^{\delta}$ is defined by :
\begin{equation}
\label{gen form approx solut gen reg meth}
x_{\alpha, reg}^{\delta} = g_\alpha^{reg}(T^*T) T^* y^\delta.
\end{equation} 
In this respect, the functions $g_\alpha^{reg}(\lambda)$ associated to Tikhonov, spectral cut-off, Showalter and the new method are defined as follows:
\begin{equation}
\label{def gen func ass to each reg meth}
g_\alpha^{tik}(\lambda) = \frac{1}{\lambda + \alpha}, \quad g_\alpha^{sc}(\lambda) = \frac{1}{\lambda} 1_{\{\lambda \geq \alpha \}}, \quad g_\alpha^{sw}(\lambda) = \frac{1-e^{-\lambda/\alpha}}{\lambda} , \quad g_\alpha^{nrm}(\lambda) = \frac{1}{\lambda + (1-\lambda^{\sqrt{\alpha}})^2}
\end{equation}
where $\lambda \in (0,a]$ with $||T^*T||_{} \leq a <1$.

Before getting into comparison of the new method to other regularization techniques, let us first point out a way of computing the regularized solution $x_{\alpha,nrm}^\delta$ of the new method.
\subsection{Computation of the regularized solution $x_{\alpha,nrm}^\delta$}
One way of computing the regularized solution $x_{\alpha,nrm}^\delta$ of the new method is by computing the singular value decomposition of operator $T$. That is to find a system $(u_k,\sigma_k,v_k)$ such that:
\begin{itemize}
\item the sequence $(u_k)_k$ forms a Hilbert basis of $X$,
\item the sequence $(v_k)_k$ forms a Hilbert basis of the closure of the range of $T$,
\item the sequence $(\sigma)_k$ is positive, decreasing and satisfies $T u_k = \sigma_k v_k$ and $T^* v_k = \sigma_k u_k$.
\end{itemize}
Given that decomposition of $T$, it is trivial to see that the operator $T^*T$ is diagonal in the Hilbert basis $(u_k)_k$. Therefore, given a function $g$ defined on the interval $(0,\sigma_1^2)$, the operator $g(T^*T)$ is nothing but the diagonal operator defined on the Hilbert basis $(u_k)_k$ by $g(T^*T) u_k = g(\sigma_k^2) u_k$. Hence given the singular value decomposition $(u_k,\sigma_k,v_k)$ of $T$, from \eqref{gen form approx solut gen reg meth} (with $reg=nrm$), the regularized solution $x_{\alpha,nrm}^\delta$ can be computed as
\begin{equation}
\label{def reg sol nrm from svd}
x_{\alpha,nrm}^\delta = \sum_{k} g_\alpha^{nrm}(\sigma_k^2)\, \langle T^* y^\delta,u_k\rangle\, u_k = \sum_{k}  \frac{\sigma_k}{\sigma_k^2 + \left(1-\sigma_k^{2\sqrt{\alpha}}\right)^2} \, \langle y^\delta,v_k \rangle \,u_k.
\end{equation}
\begin{Remark}
The above singular value decomposition of operator $T$ is only possible if $T$ is a compact operator. However, it is important to notice that the new method does not apply only to compact operator. Indeed, the new method is based on the spectral family $(E_\lambda)_\lambda$ associated to the self adjoint operator $T^*T$, and spectral family exists even for non-compact operator as pointed out in \cite[Proposition 2.14]{book_engl}. This allows for the definition of a function applied to a self-adjoint non compact operator. Of course, one might ask how we can compute the regularized solution $x_{\alpha,nrm}^\delta$ in such a case. By noticing that in practice, we always discretize Equation \eqref{main equation} into matrix formulation, we can compute the singular value decomposition of the matrix representing the discretization of operator $T$ and then apply \eqref{def reg sol nrm from svd} to compute $x_{\alpha,nrm}^\delta$.
\end{Remark}
It is important to notice that a crucial step in the computation of the regularized solution  $x_{\alpha,nrm}^\delta$ is the singular value decomposition step which should be done rigorously especially for exponentially ill-posed problems. That is why we propose a state of the art algorithm as LAPACK's \texttt{dgesvd()} routine for SVD computation (see e.g. \cite[Section 8.6]{book_matrix} for description of method). For an easy application, it is to be noted that this routine is implemented in the function \texttt{svd()} in \texttt{Matlab}. In Section \ref{section_num_illust}, we will see that even for a very ill-conditionned matrix, we can still compute the regularized solution $x_{\alpha,nrm}^\delta$ very efficiently using the function \texttt{svd()} in \texttt{Matlab}.

Above, we saw that the new approximate solution $x_{\alpha,nrm}^\delta$ is computable using the singular value decomposition of operator $T$ which might be delicate to compute. However, in some cases, there is an alternative for computing $x_{\alpha,nrm}^\delta$ when the operator $\log(T^*T)$ is explicitly known. Indeed, if the operator $\log(T^*T)$ is explicitly known, then the solution~$u: \Rb_+ \to X$ to the initial value problem:
\begin{eqnarray}
\label{ODE nrm}
\begin{cases}
u'(t) - \log(T^*T) u(t)  =  0, \quad t \in \mathbb{R}_+ \\
\qquad \qquad \qquad \quad u(0) = x, 
\end{cases}
%\quad \text{or} \quad 
%\begin{cases}
%(\log(T^*T))^{-1} u'(t) -  u(t)  =  0, \quad t \in \mathbb{R}_+ \\
%\qquad \qquad \qquad \qquad \,\,\,\,\, u(0) = x, 
%\end{cases}
\end{eqnarray}
evaluated at $t = \sqrt{\alpha}$ is nothing but $(T^*T)^{\sqrt{\alpha}} x$, i.e. $(T^*T)^{\sqrt{\alpha}} x = u(\sqrt{\alpha})$. Hence, through the resolution of the ordinary differential equation \eqref{ODE nrm}, the penalty term $\norm{\left[I - (T^*T)^{\sqrt{\alpha}} \right] x}_{}^2$ can be computed and this allows to compute the approximate solution $x_{\alpha,nrm}^\delta$.

An example of exponentially ill-posed problems for which the operator $\log(T^*T)$ is known is the backward heat equation. More precisely, let $\Omega$ be a smooth subset of $\Rb^n$ with $n \leq 3$ and $u: \Omega \times (0,\bar{t}] \to \Rb$ be the solution to the initial boundary value problem
\begin{equation}
\label{heat equation}
\begin{cases}
\qquad \frac{\partial u}{\partial t}  =  \Delta u , \quad \Omega \times (0,\bar{t}) \\
\,\, u(\cdot,0) \,\,\, =f, \quad\,\,\, \Omega  \\
\qquad \quad u = 0 \quad \text{or} \quad \frac{\partial u}{\partial \nu}  =0, \quad \text{on} \quad \partial \Omega \times (0,\bar{t}].
\end{cases}
\end{equation}
Assume we want to recover the initial temperature $f \in L^2(\Omega)$ given the final temperature $u(\cdot,\bar{t})$. By interpreting the heat equation \eqref{heat equation} as an ordinary differential equation for the function 
$ U: \left[ 0,\bar{t} \right] \to \mathcal{D}(\Delta) \subset L^2(\Omega)$, $t \to U(t) = u(\cdot,t)$, with the initial value $U(0) = f$ where 
$$
\mathcal{D}(\Delta) = H^2(\Omega)\cap H_0^1(\Omega) \quad \text{or} \quad \mathcal{D}(\Delta) = \left\lbrace f \in H^2(\Omega), \quad \frac{\partial f}{\partial \nu} = 0 \,\,\, \text{on} \,\,\, \partial \Omega \right\rbrace,
$$
we get that $U(t) = \exp{(t \Delta)} f$ for $t \in (0,\bar{t}]$, where $(\exp{(t \Delta)})_{t>0}$ is the strongly continuous semi-group generated by the unbounded self-adjoint linear operator $\Delta$.
This implies that the equation satisfied by the initial temperature $f$ is nothing but 
\begin{equation}
\label{eq sat initial temp}
\exp(\bar{t}\Delta) f = u(\cdot,\bar{t}).
\end{equation}
From \eqref{eq sat initial temp}, we deduce that $T^*T = \exp{(2\bar{t}\Delta)}$ and $\log(T^*T) = 2 \bar{t} \Delta$ and thus operator $(T^*T)^{\sqrt{\alpha}}$ can be evaluated at a function $x \in L^2(\Omega)$ as the solution to the initial value problem
\begin{equation}
\begin{cases}
u'(t) - 2 \bar{t} \Delta u(t)  =  0, \quad t \in \mathbb{R}_+ \\
\,\,\,\, \, \qquad \qquad u(0) = x, \\
u(t) \in \mathcal{D}(\Delta), \quad  \text{for} \quad t \in \mathbb{R}_+,
\end{cases}
\end{equation}
evaluated at $t = \sqrt{\alpha}$. 

In addition to the backward heat equation, there are other exponentially ill-posed problems for which $\log(T^*T)$ is known. This includes sideways heat equation 
%\footnote{For the sideways heat equation, we have $T^*T = \mathcal{F}^{-1} |\cosh\sqrt{i \omega}|^{-2} \mathcal{F}$ and $\log(T^*T) = \mathcal{F}^{-1} \log(|\cosh\sqrt{i \omega}|^{-2}) \mathcal{F}$ where $\mathcal{F}$ denotes the Fourier transform.} 
(see~\cite[Section 8.3]{reg_exp_ipp}) 
  and more generally inverse heat conduction problems (see, e.g. \cite[Section 3 \& 4]{book_murio}).
\subsection{Tikhonov versus new method}
From the variational formulation of Tikhonov and the new method, we can see that both methods differ by the penalty term. For Tikhonov method, the penalty term is $\alpha\norm{x}_{}^2$ whereas for the new method, the penalty term is $\left\| \left[I - (T^*T)^{\sqrt{\alpha}} \right] x \right\|_{}^2$ . By considering $X = L^2(\Rb^n)$ for instance, by using the Parseval identity, we see that the penalty term is equal to $\alpha \left\| \hat{x} \right\|_{L^2(\Rb^n)}$. Therefore the weight $\alpha$ equally penalizes all frequency components irrespective of the magnitude of frequencies even though instability mainly comes from high frequency components. This is actually a drawback of the Tikhonov method which may induce an unfavorable trade-off between stability and fidelity to the model (see e.g.~\cite{Alibaud}, Figure \ref{Figure num_asp}). On the contrary, for the new regularization method, high frequency components are much more regularized compared to low frequency components which are less and less regularized as the singular values increase to $1$. In this way, we expect the new method to achieve a better trade-off between stability and fidelity to the model. Moreover, for exponentially ill-posed problems, the ill-posedness is accentuated due to the magnitude of singular values, the instability introduced by high frequency components are more pronounced and we expect the new regularization method to yield better approximations of $x^\dagger$.

\subsection{Spectral cut-off versus new method}
On the contrary to Tikhonov method, both spectral cut-off and the new method treat high frequency components and low frequency components differently. However, spectral cut-off regularized high frequency components by a mere cut-off and this may be too violent in several situations. Indeed even though high frequency components induce instability, they also carry some information which should not completely left out. For instance, for mildly ill-posed problems, this truncation will be very damaging on the quality of the approximation while for exponentially ill-posed problem, this truncation will be less damaging. A smooth transition (in term of regularization) from small singular values to other singular values would be more meaningful and desirable. This is actually what is done for the new method. Another advantage of the new method compared to spectral cut-off is the variational formulation of the new method which allows to add to the problem a-priori constraint on the solution (e.g. positivity, geometrical constraints, etc...).

\subsection{Showalter versus new method}
A major difference between Showalter method and the new method is that Showalter method does not have a variational formulation. Given that, for the Showalter method, it is not clear what is actually penalized in order to stabilize the problem. Moreover it would be difficult if not impossible to add a-priori constraints on the solution. Given a data $y^\delta$, by inspecting the Showalter regularized solution which is given by $
x_\alpha^\delta = \int_{0}^{1/\alpha} e^{- s T^*T } \ud s\, T^* y^\delta
$, we see that the method introduces stability by truncating the integral
$
\int_{0}^{+\infty} e^{- s T^*T } \ud s \, T^*y^\delta = (T^*T)^{-1} T^*y^\delta
$
on the interval $(0,1/\alpha)$. On the other hand, we can see that, as the Tikhonov method, for all regularization parameter $\alpha>0$, the \text{\it generator} function $g_\alpha^{sw}$ of Showalter method is strictly decreasing on the contrary to the generator function $g_\alpha^{nrm}$ of the new method which always exhibits a maximum close to $\lambda=0$. This implies that the Showalter method cannot be seen as a smooth version of spectral cut-off which yields a smooth transition (in term of regularization) from high frequency components to low frequency components, on the contrary to the new method.
%Moreover, for $\lambda \ll \alpha$, we have $g_\alpha^{sw} \approx 1/\alpha$ and $g_\alpha^{tik} \approx 1/\alpha$, for $\lambda \approx \alpha$, we have $g_\alpha^{sw}(\lambda) \approx (1-e^{-1})/\lambda \simeq 1/(2\lambda)$ and $g_\alpha^{tik}(\lambda) \approx 1/(2\lambda)$ and finally for $\lambda \gg \alpha$, $g_\alpha^{sw}(\lambda) \approx 1/\lambda$ and $g_\alpha^{tik}(\lambda) \approx 1/\lambda$.
%All these facts can be corroborated by Figure \ref{Figure num_asp} where both functions $g_\alpha^{tik}$ and $g_\alpha^{tik}$ look very similar. Therefore, the Showalter method tends to resemble more Tikhonov method which, in $L^2$ setting, applies the same regularization to all frequency components irrespective of their magnitudes on the contrary to the new method.}
Concerning the computation of the regularized solution $x_{\alpha,sw}^\delta$ for the Showalter method, it is important to notice that $x_{\alpha,sw}^\delta$ is the solution $u_\delta: \mathbb{R}_+ \to X$ of the initial value problem:
\begin{eqnarray}
\begin{cases}
\label{ODE showalter}
u_\delta'(t) + T^*T u_\delta(t)  =  T^* y^\delta, \quad t \in \mathbb{R}_+ \\
\qquad \qquad \quad u_\delta(0) = 0, 
\end{cases}
\end{eqnarray}
evaluated at $t = 1/\alpha$, i.e. $x_{\alpha,sw}^\delta = u_\delta(1/\alpha)$.
By solving \eqref{ODE showalter} using the forward finite difference of step size $h$, we get that $u_\delta$ can be approximated as:
\begin{equation}
\label{def solution showalter by ode with forward finit diff}
u_\delta(t + h) \approx u_\delta(t) + h \left[ T^*y^\delta - T^*T u_\delta(t)  \right],  \quad \text{with} \quad u_\delta(0) = 0.
\end{equation}
\subsection{Conjugate gradient versus new method}
Unlike all the other regularization methods under consideration (Tikhonov, spectral cut-oof, Showalter and the new method), the conjugate gradient method is an iterative non-linear regularization method. The conjugate gradient method regularizes Problem \eqref{main equation} by iteratively approximating $x^\dagger$ by the minimizer $x_k$ of the functional $f(x) = || T x -y||^2$ on finite dimensional Krylov subspaces 
$$
V_k = \span\left\lbrace T^*y, (T^*T)T^*y,...,(T^*T)^{k-1}T^*y \right\rbrace,
$$
where $k\geq1$ and $k \in \mathbb{N}$. A major advantage of the conjugate gradient is the easy computation of regularized solution $x_k$ (see e.g. algorithm given in \cite[Figure 2.2]{book_kirsch}) and the fast convergence on the contrary to Landweber. However, as pointed out in \cite[Theorem 7.6]{book_engl}, the operator $R_k$ which maps the data $y$ to the regularized solution $x_k$ is not always continuous contrarily to the new method. Moreover, compared to other regularization methods, there is no a-priori rules $k(\delta)$ such that $x_{k(\delta)}^\delta$ converges to $x^\dagger$ as $\delta \to 0$ (see, e.g. \cite{ref_unst_cg}).

\begin{figure}[h!]
\begin{center}
\includegraphics[scale=0.5]{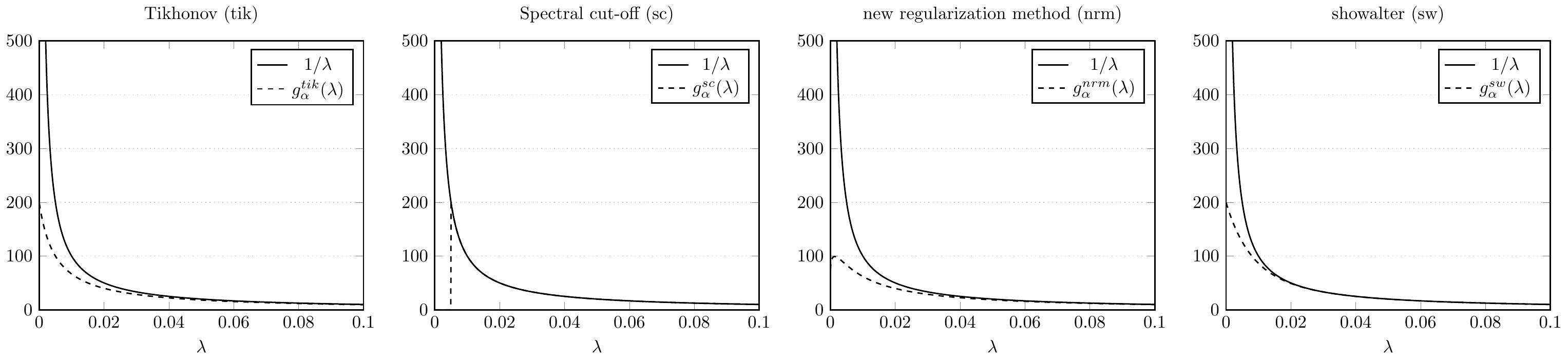} 
\end{center}
\caption{Comparison generator function $g_\alpha^{reg}$ to function $\lambda \mapsto 1/\lambda$ for the four regularization methods (reg = tik,sc,nrm,sw).}
\label{Fig_plot_gen_func_reg_mth}
\end{figure}

A comparative plot of the \text{\it generator} functions $g_\alpha^{reg}$ associated to Tikhonov, spectral cut-off, Showalter and the new method is given in Figure \ref{Fig_plot_gen_func_reg_mth}.

\begin{Remark}
On the contrary to generator functions of Tikhonov and Showalter, the generator function $g_\alpha^{nrm}$ associated to the new regularization always exhibits a maximum close to $\lambda = 0$ and the function always equals $1$ at $\lambda=0$. Indeed, it is trivial to check that both functions $g_\alpha^{tik}$ and $g_\alpha^{sw}$ are strictly decreasing for all $\alpha>0$. Hence, the function $g_\alpha^{nrm}$ is the only one which can be seen as a smooth version of the function $g_\alpha^{sc}$ associated to spectral cut-off which has a very crude transition at $\lambda = \alpha$.
\end{Remark}

%\begin{Remark}
%Given that for the moment, the computation of the regularized solution $x_{\alpha,nrm}^\delta$ for the new method relies on SVD, it is important to point out that, except for problems where the SVD is explicitly known, the new method will be difficultly applicable to large scale problems. On the contrary, the new method suits well to one dimensional problems, and hopefully two dimensional problems.
%\end{Remark}

\section{Numerical illustration}{\label{section_num_illust}}
The aim here is to compare the performance of our new regularization method (\texttt{nrm}) to the classical Tikhonov method (\texttt{tik}), spectral cut-off (\texttt{tsvd}), Showalter (\texttt{sw}) and conjugate gradient (\texttt{cg}) for some (ill-posed) test problems. We consider three test problems. The first one is a problem of image reconstruction found in \cite{shaw}. The second problem is a Fredholm integral equation of the first kind taken from \cite{baart} and the last one is an inverse heat problem. For the discretization of these problems, we use the functions \texttt{shaw()}, \texttt{baart()} and \texttt{heat()} of the \texttt{matlab} regularization tool package (see \cite{package_hansen}). 
For the \texttt{heat()} and \texttt{shaw()} test problems, the discretization is done by collocation with approximation of integrals by quadrature rules. For the \texttt{baart()} test problem, the discretization is done by Galerkin methods with orthonormal box functions as basis functions. In the \texttt{matlab} regularization tool package, each of the functions \texttt{shaw()}, \texttt{baart()} and \texttt{heat()} takes as input a discretization level $n$ representing either the number of collocations points or the number of box functions considered on the interval $[0,1]$. Given the input $n$, each function returns three outputs: a matrix $A$, a vector $x^\dagger$ and the vector $y$ obtained by discretization without noise added. In this section, we considered the following discretization level for the \texttt{shaw()}, \texttt{baart()} and \texttt{heat()} test problem respectively: $n_{shaw} = 160$, $n_{baart} = 150$ and $n_{heat} = 150$. For the simulations, we define noisy data $y_\xi = y + \xi$ where $\xi$ is a random white noise vector. In order to compute the regularized solution $x_{\alpha,nrm}^\delta$ for the new method, we compute the SVD with the function \texttt{svd()} in \texttt{Matlab} and applied \eqref{def reg sol nrm from svd}.

\begin{figure}[h!]
\begin{center}
\includegraphics[scale=1]{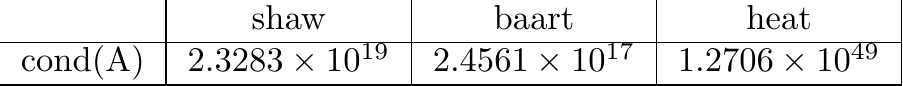} 
\end{center}
\caption{Conditioning of the matrices $A_{shaw}$, $A_{baart}$ and $A_{heat}$ for $n_{shaw} = 160$, $n_{baart} = 150$ and $n_{heat} = 150$.}
\label{Table cond}
\end{figure}

We consider a $4\%$ noise level, the noise level being defined here by the ratio of the noise to the exact data. More precisely, given a noisy data $y_\xi = y + \xi$, the noise level is defined by $\sqrt{E(||\xi||^2)}/||y||$. In order to illustrate the ill-posedness of each test problem, we give on Figure \ref{Table cond} the conditioning associated to each matrix $A_{shaw}$, $A_{baart}$, and $A_{heat}$ obtained from the discretization of each problem.

\begin{figure}[h]
   \centering
\includegraphics[scale=0.7]{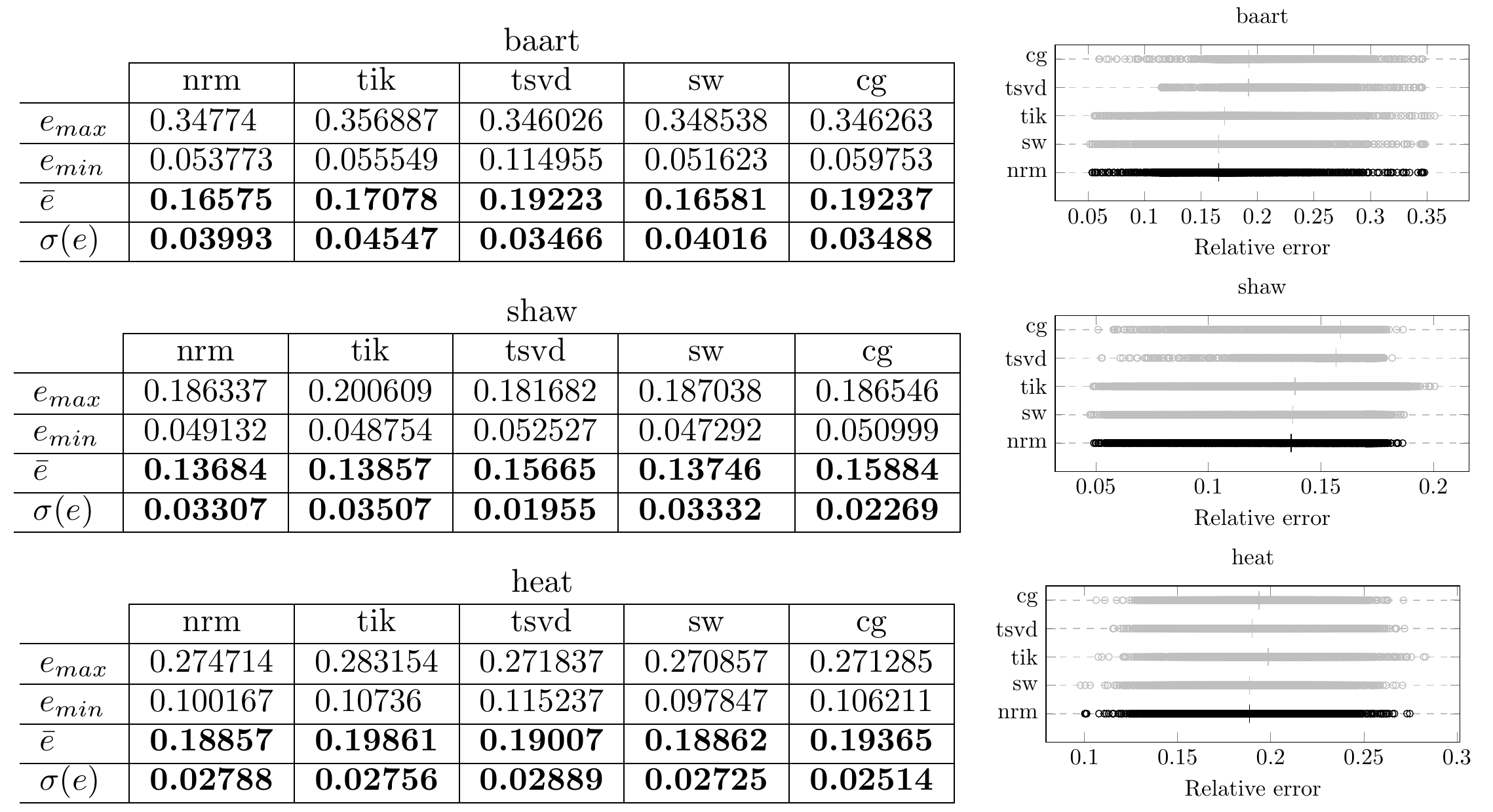}  
\caption{Summary of the Monte Carlo experiment. On the right figure, the average relative error for each method is represented by the vertical stick.}
\label{Table num_asp}
\end{figure}

 We perform a Monte Carlo experiment of 3000 replications. In each replication, we compute the best relative error for each regularization method. Next we compute the minimum, maximum, average and standard deviation errors (denoted by $e_{min}$, $e_{max}$, $\bar{e}$, $\sigma(e)$ over the 3000 replications for each schemes (\texttt{nrm} and \texttt{tik}, \texttt{tsvd}, \texttt{sw} and \texttt{cg}). Figure \ref{Table num_asp} summarizes the results of the overall simulations. 

\begin{figure}[h!]
\begin{center}
\includegraphics[scale=1]{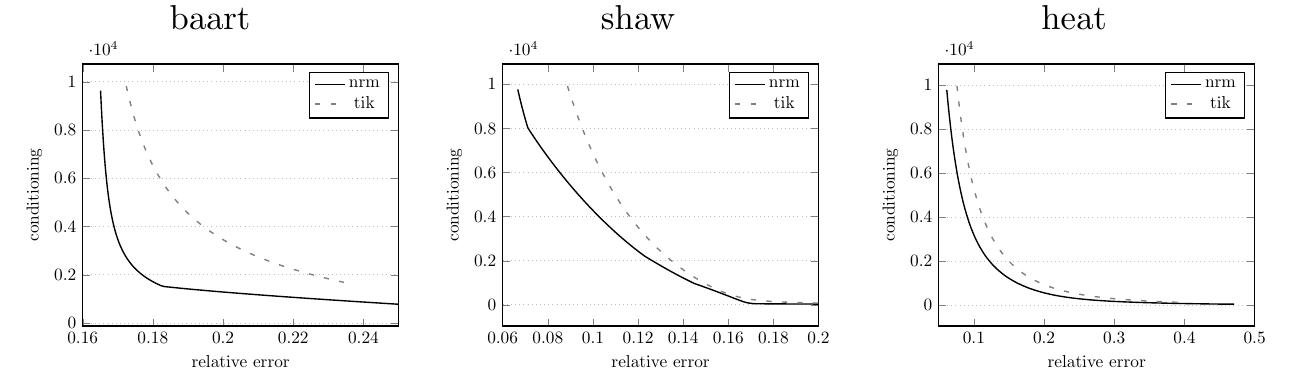}  
\end{center}
\caption{Comparison of the trade-off between stability and accuracy of the new method (\texttt{nrm}) to Tikhonov (\texttt{tik}) for the three test problems: \texttt{shaw}, \texttt{baart} and \texttt{heat}.}
\label{Figure num_asp}
\end{figure}

In order to assess and compare the trade-off between stability and fidelity to the model for Tikhonov and the new method, we plot the curve of the conditioning versus relative error.  The conditioning here is the condition number of the reconstructed operator $g_\alpha^{reg}(T^*T)$ associated to the regularization method. For instance, using the invariance of conditioning by inversion, for the new method, the conditioning corresponds to the condition number of the operator $T^*T +  \left[I - (T^*T)^{\sqrt{\alpha}} \right]^2$ while for Tikhonov method, it corresponds to the condition number of $T^*T + \alpha I$. In this respect, for two regularization methods, the best one is the one whose curve is below the other one as it achieves the same relative errors with smaller conditioning. On Figure \ref{Figure num_asp}, for each test problem, we compare the curve of conditioning versus relative error of the new method and Tikhonov method.

Notice that the first two problems (\texttt{shaw} and \texttt{baart}) are mildly ill-posed while the third problem (\texttt{heat}) is exponentially ill-posed.

\begin{flushleft}
\textbf{Comments:}
\end{flushleft}
From Figure \ref{Table num_asp} and \ref{Figure num_asp}, we can do the following comments:
\begin{itemize}

\item The new method always yields the smallest average relative errors among the five methods.

\item From Figure \ref{Table num_asp}, we can see that both spectral cut-off and conjugate gradient yield the worst average relative errors except for the \texttt{heat} test problem where their average relative errors are smaller than the one of Tikhonov.

\item For the two mildly ill-posed problem \texttt{shaw} and \texttt{baart}, Tikhonov method yields average relative errors close to the smallest one. On the contrary, for the exponentially ill-posed problem \texttt{heat}, Tikhonov method yields the worst average relative error among all the five methods.

\item For the two mildly ill-posed problems (\texttt{shaw} and \texttt{baart}),  the errors of the new method are not significantly smaller than those of Tikhonov on the contrary to the exponentially ill-posed problem (\texttt{heat}) where the new method produces smaller error than Tikhonov (about $ 5\%$ smaller). This confirms our prediction about the better performance of the new method in instance of exponentially ill-posed problems compared to Tikhonov.

\item  For all three test problems, the new method performs better than spectral cut-off as could be expected. Moreover, the gap between the error is larger for the first two  test problems which are mildly ill-posed. This also confirms the prediction about the poor performance of spectral cut-off for mildly ill-posed problems.

\item On the contrary to the two mildly ill-posed problems (\texttt{shaw} and \texttt{baart}), spectral cut-off performs better than Tikhonov  on the last test problem (\texttt{heat}), which is exponentially ill-posed. This emphasizes, especially in exponentially ill-posed problems, the drawback of Tikhonov method which regularizes all frequency component in the same way.

\item From Figure \ref{Figure num_asp}, we can see that the new method achieves a better trade-off between stability and fidelity to the model compared to the Tikhonov method. Indeed, for the three test problems the curve associated to the new method lies below the one of Tikhonov. This means that given a stability level $\kappa$ (measured in term of conditioning), the new method provided a smaller error than Tikhonov. Conversely, for a given error level $\epsilon$, the new method provides a lower conditioning of the reconstructed operator compared to Tikhonov. This also validates the prediction stated earlier.
\end{itemize}

\section{Parameter selection rules}{\label{Section_par_sel_rules}}
In this section, we are interested in the choice of the regularization parameter $\alpha$. For practical purposes, we assume that we don't know the smoothness conditions satisfied by the unknown solution~$x^\dagger$. Consequently, we are left with two types of parameter choice rules: A-posteriori rules which use information on the noise level $\delta$ and heuristic rules which depend only on the noisy data~$y^\delta$. However a huge default of a-posteriori parameter choice rules is their dependence on the noise level $\delta$ which, in practice, is hardly available or well estimated in most circumstances.
% Besides the two a-posteriori parameter choice rules defined in \eqref{morozov_like_par_sel_rule} and \eqref{morozov_simple_par_sel_rule}, an interesting converging a-posteriori parameter choice rule that could suit for the new method is the one proposed in \cite{geom_lin_ipp} by Math\'e and Pereverzev which is based on the Lepskii approach \cite{paper_lepskii}. The first step of the method consists in discretizing the regularization parameter as follows:
%\begin{equation*}
%\alpha_n = \alpha_0 q^n, \quad \text{with} \quad \alpha_0 > \delta^2 \quad, q>1, \quad \text{and} \quad n=1,2,...N
%\end{equation*}
%where $N$ satisfies $\alpha_{N-1} \leq ||T^*T|| \leq \alpha_N$.
%Next the regularization parameter  $\bar{\alpha}$ is chosen as 
%\begin{equation}
%\label{reg par rule mathe and pereverzev}
%\bar{\alpha} := \max \left\lbrace \alpha_i, \quad ||x_{\alpha_i}^\delta - x_{\alpha_{i-1}}^\delta || \leq 4 C_\gamma \frac{\delta}{\sqrt{\alpha_{i-1}}}  \right\rbrace,
%\end{equation}
%where $ C_\gamma = \max \{\gamma_*, \bar{\gamma} \}$.  Here $\gamma_*$ is the constant on the right hand side in \eqref{bound fct prop err} (which can be computed) and $\bar{\gamma}$ is the constant on the right hand side in the qualification inequality, e.g. for logarithmic source conditions, $\bar{\gamma}$ is the constant in the right hand side of \eqref{qualif_new_meth}. However, the unavailability of the constant $\bar{\gamma}$ for the new method makes this rule not practical for the new method.
 In \cite{note_morozov}, it is shown how an underestimation or overestimation of the noise level $\delta$ may induce serious computation issues for the Morozov principle. Moreover, in \cite{reg_large_scale_prob}, it is illustrated how heuristic rules may outperform sophisticated a-posteriori rules. Given those reasons, we turn to heuristic (or data driven) selection rules. We recall that, due to Bakushinskii v\'eto \cite{non_conv_heur_rule}, such rules are not convergent. But still, as mentioned earlier, heuristic rules may yields better approximations compared to sophisticated a-posteriori rules (see e.g. \cite{reg_large_scale_prob}) and this is not surprising as the Bakushinskii result is based on worst case scenario.

We applied five noise-free parameter choice rules to the new method and the four regularization methods on the three test problems defined in Section \ref{section_num_illust}: the generalized cross validation (GCV), the discrete quasi-optimality rule (DQO), two heuristic rules (H1 and H2) and a variant of the L-curve method (LCV) each described in \cite[Section 4.5]{book_engl}. Roughly speaking, the parameter $\alpha$ chosen by each of those selection rules is as follows:
\begin{itemize}
\item The GCV rule consists in choosing $\hat{\alpha}$ as 
$$
\hat{\alpha} = \argmin_{\alpha} \frac{|| T x_\alpha^\delta - y^\delta||}{\tr(r_\alpha(T^*T))},
$$
where $r_\alpha$ is the \textit{residual} function associated to the regularization method under consideration. For the new method, $r_\alpha$ is defined in \eqref{def_res_func_new_method}.
\item The DQO method consists in discretizing the regularization parameter $\alpha$ as
$$
\alpha_n = \alpha_0 q^n, \quad \alpha_0 \in (0,||T^*T||], \quad \text{and} \quad 0< q <1.
$$
Next, the parameter $\hat{\alpha}$ is chosen as 
\begin{equation}
\label{DQO rule}
\hat{\alpha} = \alpha_{\hat{n}} \quad  \text{with} \quad \hat{n} = \argmin_{n \in \Nb}  ||x_{\alpha_{n+1}}^\delta  - x_{\alpha_n}^\delta||.
\end{equation}
Recall that this rule defined by \eqref{DQO rule} is actually one of the variants of the continuous quasi-optimality rule defined by 
$$
\hat{\alpha} = \argmin_\alpha || \alpha \frac{\partial x_{\alpha}^\delta}{\partial \alpha} ||.
$$
\item The third rule H1 taken in \cite[Section 4.5]{book_engl} consists in choosing the parameter $\hat{\alpha}$ as
\begin{equation}
\label{H1 rule}
\hat{\alpha} = \argmin_{\alpha} \frac{1}{\sqrt{\alpha}} || T x_\alpha^\delta - y^\delta||.
\end{equation}
\item The fourth rule H2 which is a variant of the third rule H1 consists in choosing the parameter $\hat{\alpha}$ as
\begin{equation}
\label{H2 rule}
\hat{\alpha} = \argmin_{\alpha} \frac{1}{\alpha} || T^*(T x_\alpha^\delta - y^\delta) ||.
\end{equation}
\item The variant of the L-curve (LCV) method considered here (see \cite[Proposition 4.37]{book_engl}) consists in choosing the regularization parameter $\hat{\alpha}$ as
\begin{equation}
\label{LCV rule}
\hat{\alpha} = \argmin_\alpha ||x_\alpha^\delta||\,||T x_\alpha^\delta - y^\delta||.
\end{equation}
Recall that this rule actually tries to locate the parameter $\hat{\alpha}$ corresponding to the corner of the L-curve plot $||x_\alpha^\delta||$ versus $||T x_\alpha^\delta - y^\delta||$ in a log-log scale. For more details about the L-curve method, see e.g. \cite{l_curve_morozov,l_curve_1,l_curve_2}.
\end{itemize}
For a comprehensive discussion of the above heuristic rules and conditions under which convergence is established, see \cite{gcv,asymptotic_optimality_gcv,gcv1} for GCV, \cite{recent_result_quasi_opt,quasi_opt_1,quasi_opt2,leonov}) for Quasi-optimality and \cite[Section 4.5]{book_engl} for the rules H1, H2 and LCV. 

\begin{figure}[h!]
\begin{center}
%\begin{flushleft}
\includegraphics[scale=0.4]{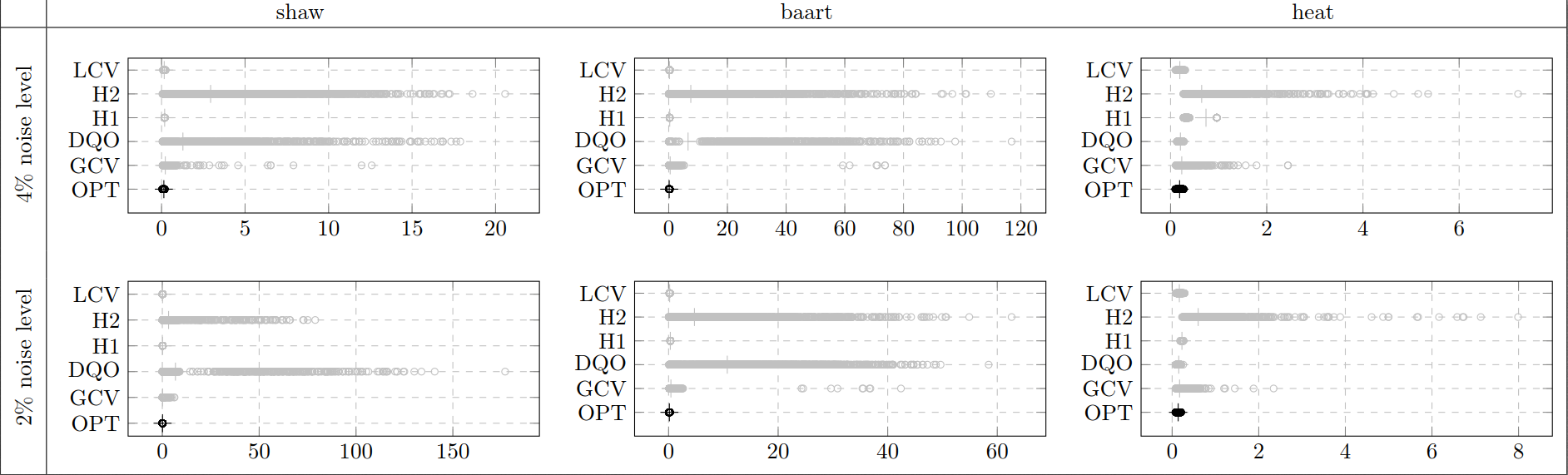} 
%\end{flushleft}
\end{center}
\caption{Comparison of the relative error obtained by each selection rules (GCV, DQO, H1, H2 and LCV) for the two noise levels with the new method for the three tests problems \texttt{shaw}, \texttt{baart} and \texttt{heat}. On each plot, the x-axis corresponds to relative error and the vertical stick indicates the average relative error.}
\label{Figure_sum_heur_rule}
\end{figure}

For assessing the performance of each selection rule, we perform a Monte Carlo experiment of 3000 replications. For each replication, each test problem (\texttt{baart}, \texttt{shaw}, \texttt{heat}), and each regularization method (\texttt{nrm}, \texttt{tik}, \texttt{tsvd}, \texttt{sw} and \texttt{cg}), we compute the optimal regularization parameter $\alpha_{OPT}$, the one chosen by each selection rule $(\alpha_{GCV}, \alpha_{DQO}, \alpha_{H1}, \alpha_{H2},\alpha_{LCV})$. We also compute the corresponding relative errors: 
$$
\frac{||x^\dagger - x_{\alpha_{OPT}}^\delta||}{||x^\dagger||}, \,\,\, \frac{||x^\dagger - x_{\alpha_{GCV}}^\delta||}{||x^\dagger||}, \,\,\, \frac{||x^\dagger - x_{\alpha_{DQO}}^\delta||}{||x^\dagger||}, \,\,\, \frac{||x^\dagger - x_{\alpha_{H1}}^\delta||}{||x^\dagger||}, \,\,\, \frac{||x^\dagger - x_{\alpha_{H2}}^\delta||}{||x^\dagger||}, \,\,\, \text{and} \,\,\, \frac{||x^\dagger - x_{\alpha_{LCV}}^\delta||}{||x^\dagger||}.
$$
In order to analyse the convergence behavior of the selection rules, we consider two noise levels: $2\%$ and $4\%$. The results are shown in Figure \ref{Figure_sum_heur_rule} and  Tables \ref{Table_sum_heur_rule_1} to \ref{Figure_lcv_heat}.

\begin{table}[h!]
\begin{center}
\includegraphics[scale=0.8]{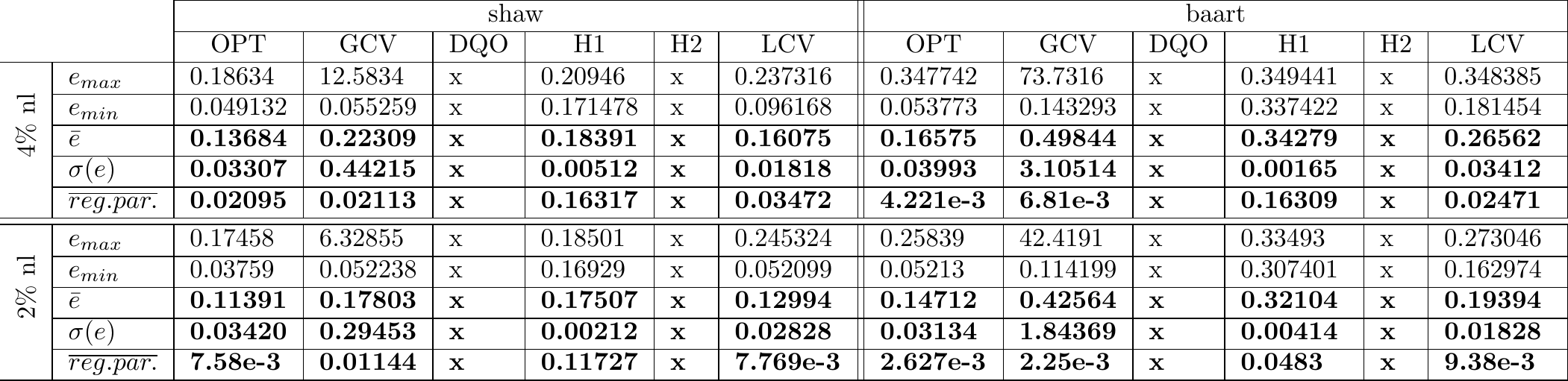} 
\end{center}
\caption{Summary of the Monte carlo experiment with the five heuristic rules GCV, DQO, H1, H2 and LCV applied to the new method for the test problems \texttt{shaw} and \texttt{baart}. The x indicates columns where the average relative error is greater than $1$.}
\label{Table_sum_heur_rule_1}
\end{table}

\begin{table}[]
\begin{center}
\includegraphics[scale=0.9]{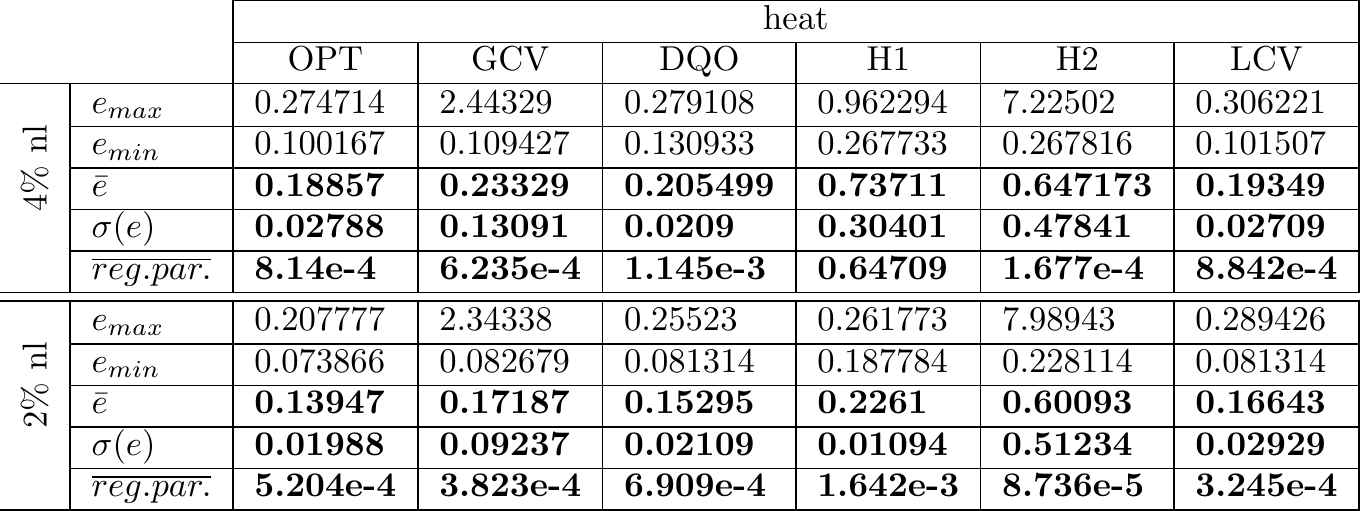}  
\end{center}
\caption{Summary of the Monte carlo experiment with the five heuristic rules GCV, DQO, H1, H2 and LCV applied to the new method for the test problem \texttt{heat}.}
\label{Table_sum_heur_rule_2}
\end{table}

\begin{table}[h!]
%\begin{center}\\
\begin{flushleft}
\includegraphics[scale=0.8]{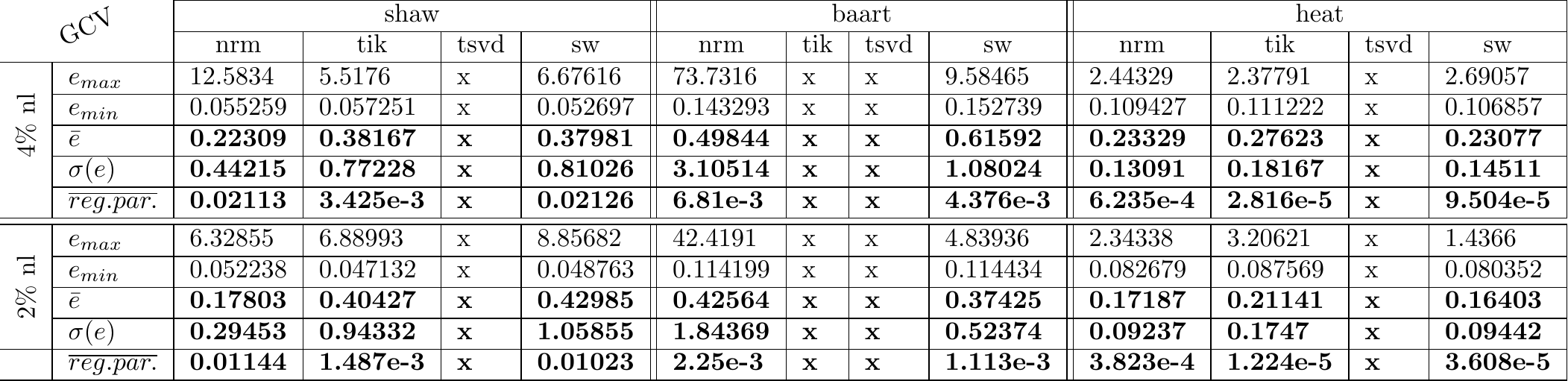} 
\end{flushleft}
%\end{center}
\caption{Summary of the Monte Carlo experiment with GCV rule applied to nrm,tik,tsvd and sw for the two noise levels on the three tests problems \texttt{shaw}, \texttt{baart} and  \texttt{heat}.The x indicates columns where the average relative error is greater than $1$.}
\label{Figure_gcv_baart_shaw}
\end{table}

\begin{table}[h!]
%\begin{center}
\begin{flushleft}
\includegraphics[scale=0.70]{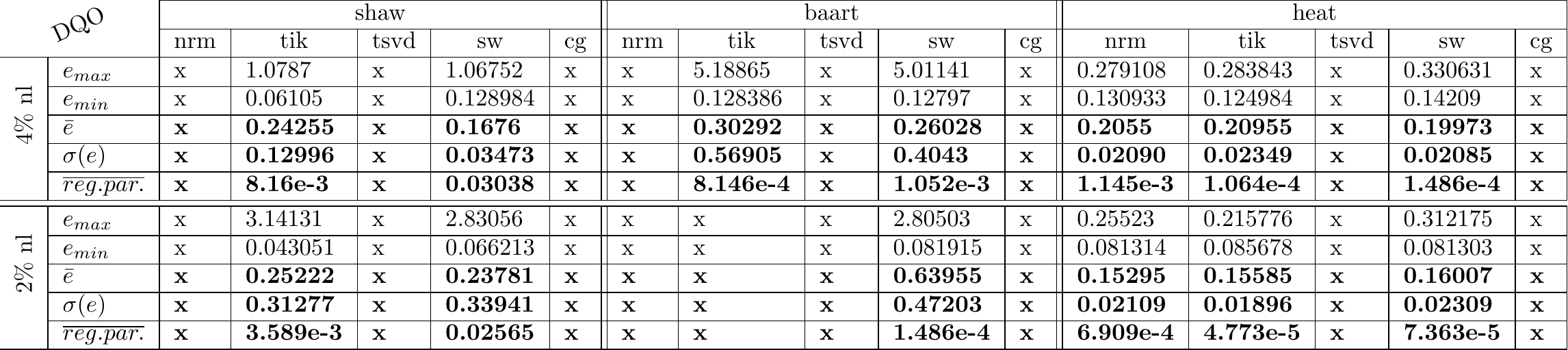} 
\end{flushleft}
%\end{center}
\caption{Summary of the Monte Carlo experiment with DQO rule applied to nrm,tik,tsvd,sw and cg for the two noise levels on the three tests problems \texttt{shaw}, \texttt{baart} and  \texttt{heat}.The x indicates columns where the average relative error is greater than $1$.}
\label{Figure_dqo_baart_shaw}
\end{table}

\begin{table}[h!]
\begin{center}
\includegraphics[scale=0.8]{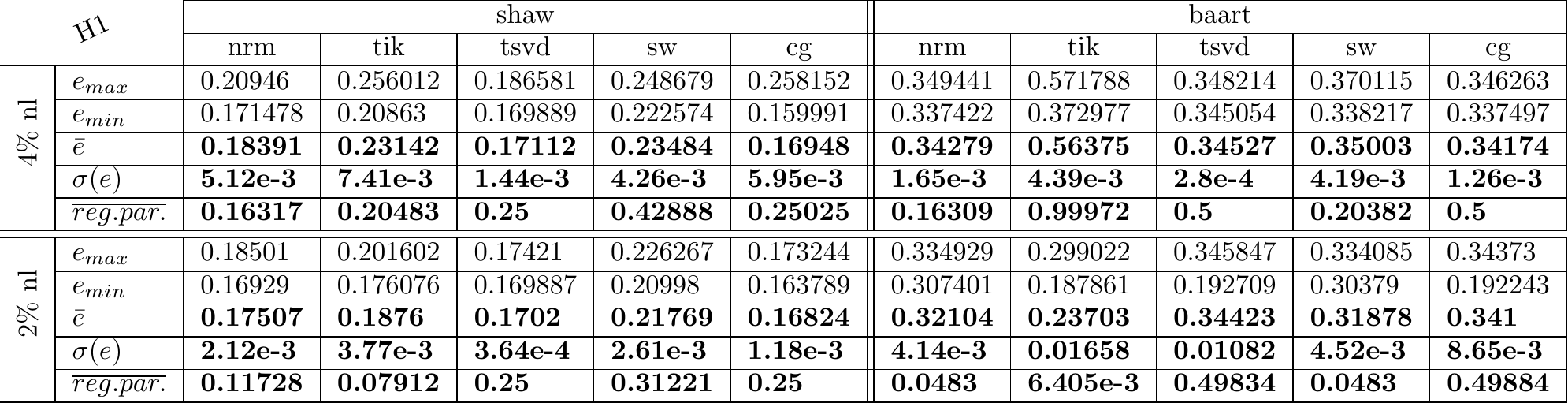}
\end{center}
\caption{Summary of the Monte Carlo experiment with rule H1 applied to nrm,tik,tsvd,sw and cg for the two noise levels on the two tests problems \texttt{shaw} and \texttt{baart}}
\label{Figure_h1_baart_shaw}
\end{table}

\begin{table}[h!]
%\begin{center}
\begin{flushleft}
\includegraphics[scale=0.7]{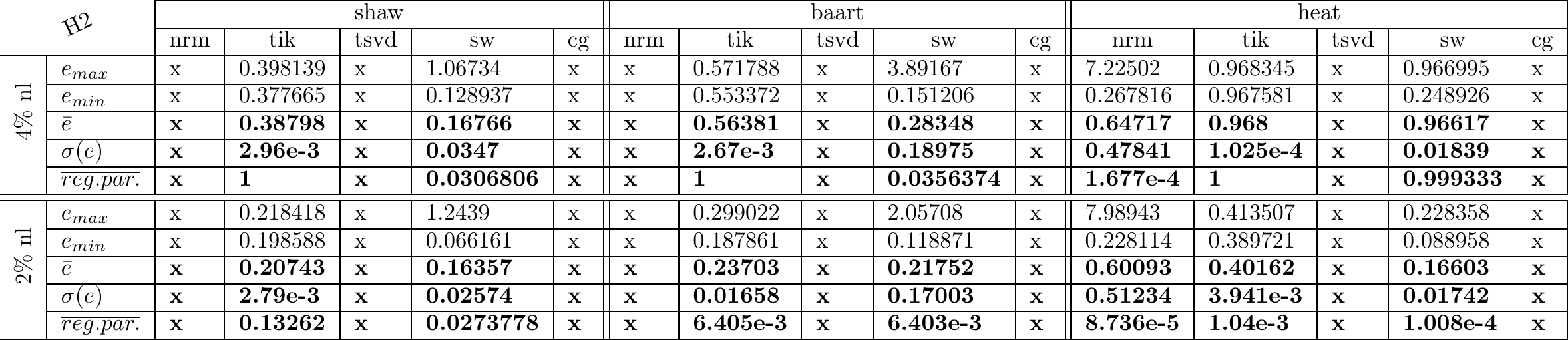}
\end{flushleft}
%\end{center}
\caption{Summary of the Monte Carlo experiment with rule H2 applied to nrm,tik,tsvd,sw and cg for the two noise levels on the three tests problems \texttt{shaw}, \texttt{baart} and  \texttt{heat}.The x indicates columns where the average relative error is greater than $1$.}
\label{Figure_h2_baart_shaw}
\end{table}

\begin{table}[h!]
\begin{center}
\includegraphics[scale=1]{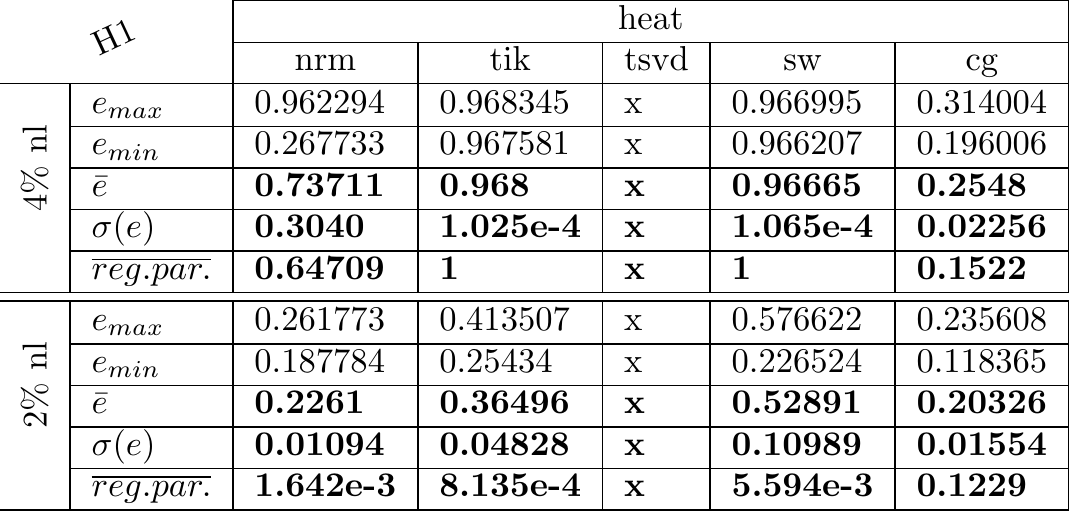} 
\end{center}
\caption{Summary of the Monte Carlo experiment with rule H1 applied to nrm,tik,tsvd,sw and cg for the two noise levels on the test problem \texttt{heat}.The x indicates columns where the average relative error is greater than $1$.}
\label{Figure_h1_heat}
\end{table}

\begin{table}[h!]
\begin{center}
\includegraphics[scale=0.8]{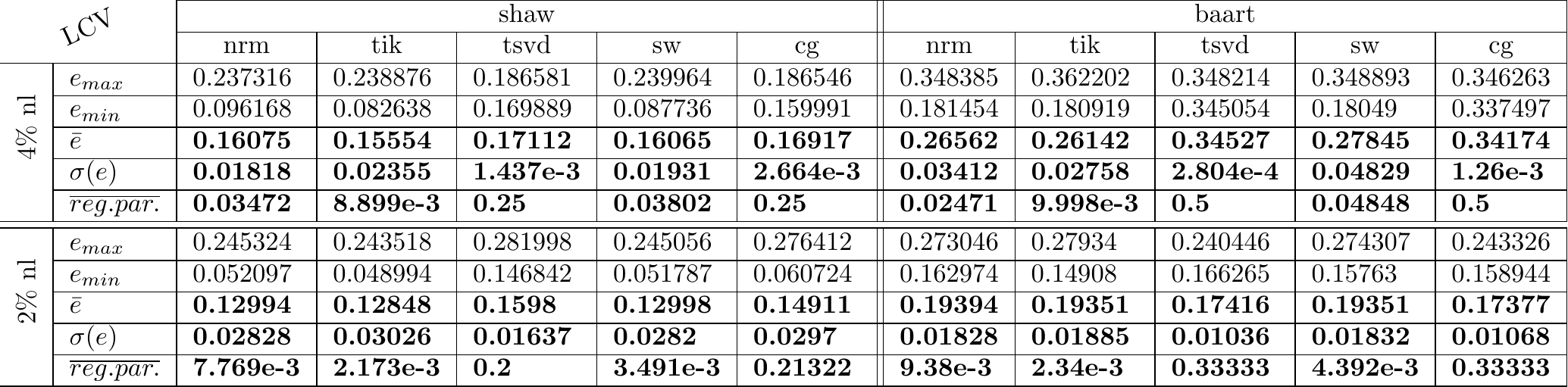}
\end{center}
\caption{Summary of the Monte Carlo experiment with LCV rule applied to nrm,tik,tsvd,sw and cg for the two noise levels on the two tests problems \texttt{shaw} and \texttt{baart}}
\label{Figure_lcv_baart_shaw}
\end{table}

\begin{table}[h!]
\begin{center}
\includegraphics[scale=1]{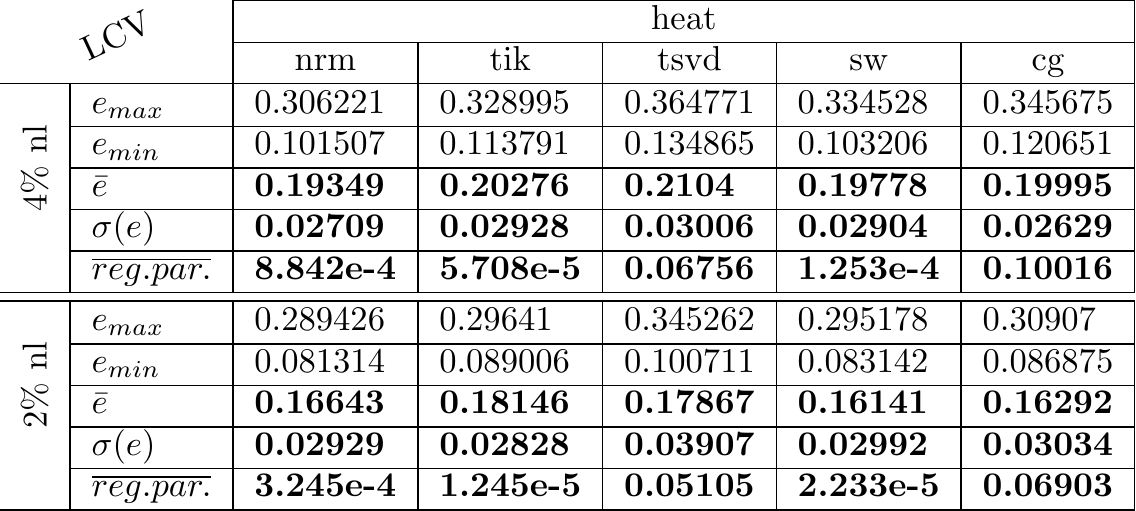} 
\end{center}
\caption{Summary of the Monte Carlo experiment with LCV rule applied to nrm,tik,tsvd,sw and cg for the two noise levels on the test problem \texttt{heat}.}
\label{Figure_lcv_heat}
\end{table}

From Tables \ref{Table_sum_heur_rule_1}, \ref{Table_sum_heur_rule_2} and Figure \ref{Figure_sum_heur_rule}, we can see the following concerning the new regularization method: 
\begin{itemize}
\item For the exponentially ill-posed problem \texttt{heat}, from Table \ref{Table_sum_heur_rule_2} and the last column of Figure \ref{Figure_sum_heur_rule},  we can see that the discrete quasi-optimality rule and the variant of the L-curve are very efficient parameter choice rules for the new method. Indeed both the average relative errors and the average regularization parameters produced by the DQO and LCV rules are very near the optimal ones and decrease as the noise level decreases. Moreover, by looking at the standard deviation of the relative error $\sigma(e)$, we see that those rules are very stable with respect to variations of the error term in $y$. Next, the GCV rule exhibit good average relative error, however the GCV is not stable with respect to the noise in $y$ and this is shown by the spreading of dots along the x-axis or the corresponding large standard deviation $\sigma(e)$. Finally, the rule H2 is unstable and produces large relative errors norm whereas the rule H1 is more stable but do not yield satisfactory errors.

\item For the mildly ill-posed test problems \texttt{shaw} and \texttt{baart}, the best heuristic rule for the new method is the variant of the L-curve method. Indeed, from Table \ref{Table_sum_heur_rule_1} and two first columns of Figure \ref{Figure_sum_heur_rule}, we notice that the relative errors produced by the LCV rule are near the optimal ones. Moreover, the LCV rule is very stable with respect to the noise in $y$ and both the relatives errors and the regularization parameters decrease as the noise level decreases. The second best rule is rule H1 which is also stable and convergent but produces relative errors larger than the one of LCV rule. Finally the rules DQO, GCV and H2 are unstable and produce large relative error norm.

\end{itemize}
From Tables \ref{Figure_gcv_baart_shaw} to \ref{Figure_lcv_heat}, we apply the five selection rules GCV, DQO, H1, H2 and LCV to each regularization method. Obviously the GCV rule cannot be applied to conjugate gradient method due to its non-linear character. Although, the DQO is originally designed for continuous regularization methods, notice that the rule defined in \eqref{DQO rule} can be applied to regularization  methods with discrete regularization parameter such as truncated singular value decomposition and conjugate gradient. Indeed, we can applied the DQO rule to \texttt{tsvd} and \texttt{cg} by replacing $x_{\alpha_{n}}^\delta$ by $x_k^\delta$ in \eqref{DQO rule}.
Similarly the rules H1 and H2 originally designed for continuous regularization methods may be applicable to discrete regularization by defining the regularization parameter $\alpha$ as the inverse of the discrete parameter $k$. Following that idea, we applied the rules H1 and H2 to \texttt{tsvd} and \texttt{cg} by replacing $\alpha$ by $1/k$ in \eqref{H1 rule} and \eqref{H2 rule}.

From Tables \ref{Figure_gcv_baart_shaw} to \ref{Figure_lcv_heat}, we can do the following comments:
\begin{itemize}
\item The variant of the L-curve method defined through \eqref{LCV rule} is a very efficient heuristic parameter choice rule for each considered regularization method. Indeed, from Tables \ref{Figure_lcv_baart_shaw} and \ref{Figure_lcv_heat}, by looking at the standard deviation $\sigma(e)$ of the relative error, we see that the LCV rule is stable for each regularization method, each test problem and each noise level. Next, the rule exhibits a convergent behavior for each test problem and each regularization method since the average relative error $\bar{e}$ and the average regularization parameter $\overline{reg. par.}$ decrease as the noise level decreases. Finally from Tables  \ref{Figure_gcv_baart_shaw} to \ref{Figure_lcv_heat}, we find that the LCV rule always yields the smallest average relative error $\bar{e}$ among all the heuristic rules considered except in 4 cases (out of 30 cases in total) : \texttt{baart} test problem with $4\%$ noise level for Showalter method and \texttt{heat} test problem with $2\%$ noise level for the new method, Tikhonov and Showalter method. Notice that in each of those four cases, LCV rule yields the second best average relative error $\bar{e}$ after the DQO rule. 
\item For the exponentially ill-posed test problem \texttt{heat}, Table \ref{table best heur rule for heat test pb} summarizes the best heuristic rules for each regularization method:
\begin{table}[h!]
\begin{center}
\includegraphics[scale=1]{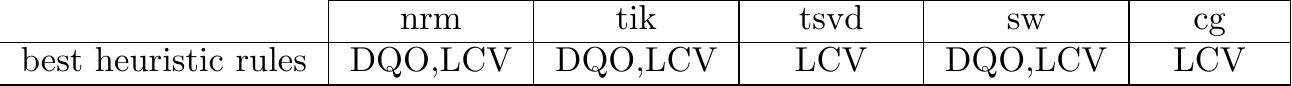} 
\end{center}
\caption{Summary best heuristic rules for each regularization method for the exponentially ill-posed test problem \texttt{heat}.}
\label{table best heur rule for heat test pb}
\end{table}
\item For the mildly ill-posed test problems \texttt{shaw} and \texttt{baart}, the best heuristic rule is always the LCV rule. For the new method, Tikhonov, truncated singular value decomposition and conjugate gradient, the LCV rule is followed by rule H1 whereas for the Showalter method, the LCV rule is followed by rule H2.
\item For the exponentially ill-posed test problem \texttt{heat}, by comparing the five regularization methods combined each with its best heuristic selection rule among GCV, DQO, H1, H2 and LCV, we see that the new method equipped with the DQO rule (resp. the LCV rule) for $4\%$ noise level (resp. for $2\%$ noise level) yields the smallest average relative error $\bar{e}$ (about $2\%$ smaller than the second best average relative error). For $4\%$ noise level, the second smallest average relative error is achieved by Showalter method equipped with LCV rule whereas for $2\%$ noise level,  the second smallest average relative error is achieved by Tikhonov method equipped with DQO rule. 
\item For the two mildly ill-posed problems \texttt{shaw} and \texttt{baart}, by comparing the five regularization methods combined each with its best heuristic selection rule among GCV, DQO, H1, H2 and LCV, we notice there is no regularization method which always yields the smallest average relative error. For the \texttt{shaw} test problem, Tikhonov method with LCV rule yields the smallest average relative error $\bar{e}$. For the \texttt{baart} test problem, for $4\%$ noise level, the smallest average relative error is obtained by the Showalter method equipped with the DQO rule. However, for this test problem, the DQO rule is not converging for the Showalter method as the average relative error $\bar{e}$ increases from $0.26028$ to $0.63955$ as the noise level decreases from $4\%$ to $2\%$. If we discard Showalter with DQO rule, then for $4\%$ noise level, the smallest average relative error are obtained by Tikhonov method equipped with LCV rule while for $2\%$ noise level, the smallest average relative errors are obtained from conjugate gradient method equipped with LCV rule.
\end{itemize}
\begin{Remark}
From Tables \ref{Table_sum_heur_rule_1} to \ref{Figure_lcv_heat}, we see that, the heuristic parameter choice rule LCV yields very satisfactory results for each considered regularization method. This reinforces the idea that the Bakushinskii v\'eto \cite{non_conv_heur_rule} should not be seen as a limitation of heuristic parameter choice rule but rather as a safeguard to be taken into account.
\end{Remark}

 In summary,  we see that for the exponentially ill-posed test problem \texttt{heat}, the new regularization method always yields the smallest average relative error among the five considered regularization methods even when we consider heuristic parameter choice rules. Hence in practical situation of exponentially ill-posed problems, we expect the new method to perform better than the other regularization methods (Tikhonov, truncated singular value decomposition, Showalter method and conjugate gradient).

%%%%%%%%%%%%%%%%%%%%%%%%%%%%%%%%%%%%%%%%%%%%%%%%%%%%%%

\section{Conclusion}
In this paper, we presented a new regularization method which is particularly suitable for linear exponentially ill-posed problems. We study convergence analysis of the new method and we provided order optimal convergence rates under logarithmic source conditions which has a natural interpretation in term of Sobolev spaces for exponentially ill-posed problems. For a general source conditions expressed via index functions, we only provided quasi order optimal rates. 
From the simulatins performed, we saw that the new method performs better than  Tikhonov method, spectral cut-off, Showalter and conjugate gradient for the considered exponentially ill-posed problem, even with heuristic parameter choice rules. For the two  mildly ill-posed problems treated, we saw that the new method actually yields results quite similar to those of Tikhonov and Showalter methods. 
The results of Section \ref{Section_par_sel_rules}, where we applied five \text{\it error-free} selection rules to the five regularization methods suggest that the variant of the L-curve method defined in \eqref{LCV rule} and the discrete quasi-optimality rule defined in \eqref{DQO rule} are very efficient parameter choice rules for the new method in the context of exponentially ill-posed problem. In the context of mildly ill-posed problems, the results of experiments suggest that the LCV rule described in Section \ref{Section_par_sel_rules} is preferable.

Interesting perspectives would be 
%to seek alternatives for computing the regularized solution of the new method without having to compute SVD of operator $T$ and 
a theoretical analysis of the LCV and DQO rules for the new regularization method in the framework of exponentially ill-posed problems in order to shed light on their good performances.
\newline

\textbf{Acknowledgements:}
The author would like to thank Pierre Mar\'echal, Anne Vanhems for their helpful comments, readings and remarks.
% and particularly the referees for their reports which helped me to significantly improve the paper.  

\section{Appendix}

\textbf{Proof of Proposition \ref{prop 2}.}
 Let us state the following standard inequality that we will use in the sequel:
\begin{equation}
\label{prop exp}
\forall \, t \geq 0,\quad \exp(-t) \leq \frac{1}{1+t}.
\end{equation}
Using \eqref{prop exp} applied with $t = - \sqrt{\alpha} \ln(\lambda) \geq 0 $, we get
\begin{equation}
\label{part demo 1}
1-\exp(\sqrt{\alpha} \ln(\lambda)) \geq 1- \frac{1}{1 - \sqrt{\alpha} \ln(\lambda)} = \frac{-\sqrt{\alpha} \ln(\lambda)}{1 - \sqrt{\alpha} \ln(\lambda)} = \sqrt{\alpha} \frac{|\ln(\lambda)|}{1 + \sqrt{\alpha} |\ln(\lambda)|}.
\end{equation}
But since $\alpha <1$, $1 + \sqrt{\alpha} |\ln(\lambda)| < 1 + |\ln(\lambda)|$. Furthermore For all $\lambda \leq a <1$, by the monotonicity of the function $t \To |\ln(t)|/(1+|\ln(t)|) = -\ln(t)/(1-\ln(t))$ on $(0,1)$, we get that 
$$
\frac{|\ln(t)|}{1+|\ln(t)|} \geq \frac{|\ln(a)|}{1+|\ln(a)|} \quad \forall \, t \in (0,a).
$$
By applying the above inequality to \eqref{part demo 1} and taking the square, we get:
$$
 \forall \lambda \in (0,a), \quad (1-\lambda^{\sqrt{\alpha}})^2 \geq M \alpha  \quad \text{with} \quad M = \left(\frac{|\ln(a)|}{1+|\ln(a)|}\right)^2.
$$
Whence the following inequality:
\begin{equation}
\label{bound gen func nrm}
\frac{1}{\lambda + (1-\lambda^{\sqrt{\alpha}})^2} \leq \frac{1}{\lambda + M\alpha},
\end{equation}
which implies that
\begin{equation}
\label{tttt}
\sqrt{\lambda} g_\alpha(\lambda) \leq \frac{\lambda^{1/2}}{\lambda + M \alpha}.
\end{equation}
It is rather straightforward to prove that the supremum over $\lambda \in (0,1)$ of the right hand side of \eqref{tttt} is of order $\alpha^{-1/2}$ from which we deduce that
\begin{equation}
\label{sharp bound nu_beta}
\sup_{\lambda \in (0,a]} \sqrt{\lambda} g_\alpha(\lambda)  = \mathcal{O}\left( \frac{1}{\sqrt{\alpha}}\right)
\end{equation}
\QED
\newline

\textbf{Proof of Lemma \ref{Lemma bound resid_func}.}
Let $\lambda \in (0,1)$. On the one hand, by applying the estimate $ (1-\exp(t)) \geq -t/(1-t) $ which holds for all $t<0$ to $t = \sqrt{\alpha}\ln(\lambda)$ and by taking squares, we have:
\begin{equation}
\label{low_bound}
(1-\lambda^{\sqrt{\alpha}} )^2 \geq \frac{\alpha |\ln(\lambda)|^2}{(1+\sqrt{\alpha} |\ln(\lambda)|)^2}.
\end{equation}
On the other hand, using the estimate $t^2 \geq (1-\exp(t))^2$ valid for all $t<0$ to $t = \sqrt{\alpha}\ln(\lambda)$, we get
\begin{equation}
\label{up_bound}
 (1-\lambda^{\sqrt{\alpha}} )^2 \leq\alpha |\ln(\lambda)|^2
 \end{equation}
Now, for $\alpha \leq \lambda <1$, $|\ln(\alpha)| \geq |\ln(\lambda)|$ which implies that $\sqrt{\alpha}|\ln(\lambda)| \leq \sqrt{\alpha}|\ln(\alpha)|$. Using the estimate $t^\mu \ln (1/t) \leq \mu$ which is true for all $t$ in $(0,1)$ and every positive $\mu$ to $t = \lambda$ and $\mu = 1/2$, we deduce that
\begin{equation}
\label{ee}
1 + \sqrt{\alpha}|\ln(\alpha)|   \leq 3/2.
\end{equation}
So, from \eqref{low_bound} and \eqref{ee}, we deduce that
\begin{equation}
\label{best_low_bound}
(1-\lambda^{\sqrt{\alpha}} )^2 \geq \frac{4}{9} \alpha |\ln(\lambda)|^2
\end{equation}
which implies that
\begin{equation}
\label{ff}
r_\alpha(\lambda) \leq \frac{(1-\lambda^{\sqrt{\alpha}} )^2}{\lambda + (4/9)\alpha |\ln(\lambda)|^2}
\end{equation}
Finally, applying \eqref{up_bound} and the fact that $\lambda \geq (4/9)\lambda$ to \eqref{ff} yields \eqref{bound_res_func_new_reg_meth}.
\QED
\newline

\textbf{Proof of Lemma \ref{Lemma 2}.}
(i) It is straightforward to check that \eqref{der_func_chi} is indeed the derivative of the function $\Psi_{p,\alpha}$. 

(ii) First notice that 
$\lim_{\lambda \to 0} h(\lambda) = + \infty $, hence, it suffices to find a $\bar{\lambda}$ such that $h(\bar{\lambda}) <0$ to deduce the existence of a root of the function $h$ on $(0,\bar{\lambda}]$. 
If $p<2$, then $h(1) <0$. If $p=2$, then $h(\lambda) = |\ln(\lambda)|(2 \alpha |\ln(\lambda)| - \lambda)$. Thus, $h(\lambda)<0$ for $\lambda$ close to $1$ but smaller than $1$. If $p>2$, then $\lim_{\alpha \to 0} h(\lambda) = \lambda (p-2 + \ln(\lambda))<0$ for all $\lambda<\exp(2-p)$.

Now let us show that for every $\lambda(p,\alpha)$ which vanishes $h$, \eqref{estimate root of h} holds.
\begin{equation}
\label{def alpha}
h(\lambda) = 0 \quad \Longrightarrow \quad \alpha = \lambda |\ln(\lambda)|^{-1} \left( \frac{2-p + |\ln(\lambda) |}{ p |\ln(\lambda)|} \right)
\end{equation}
by monotonicity of the function $t \to (2-p +t)/(pt)$ (irrespective of the sign of $2-p$) and $t \to  |\ln(\lambda)|$, we get that the function $l(\lambda) = \frac{2-p + |\ln(\lambda) |}{ p |\ln(\lambda)|}$ is monotonic.
If $p<2$, the function $l$ is increasing and we then get that, for all $\lambda \in (0,c]$ with $c <1$,
\begin{equation}
\label{bound 1 }
\frac{1}{p} \leq l(\lambda) \leq l(c).
\end{equation}
On the other hand, if $p \geq 2$, the function $l$ is decreasing and for $\lambda \in (0,c]$ with $c < \exp{(2-p)}$, we get
\begin{equation}
\label{bound 2 }
l(c) \leq l(\lambda) \leq 1/p.
\end{equation}
From \eqref{def alpha}, \eqref{bound 1 } and \eqref{bound 2 }, we deduce that 
\begin{equation}
\label{equiv def alpha}
h(\lambda) = 0 \quad \Longrightarrow \quad \alpha \sim \lambda |\ln(\lambda)|^{-1}.
\end{equation}
From \cite[Lemma 3.3]{optimality_under_gen_sour_cond}, we get that 
$$
\alpha \sim \lambda |\ln(\lambda)|^{-1} \Rightarrow \lambda \sim \alpha | \ln(\alpha)| (1 + \smallo{1} ) \quad \text{for} \quad \alpha \to 0.
$$
This shows that the maximizers $\lambda(p,\alpha)$ of the function $\Psi{p,\alpha}$ satisfies \eqref{estimate root of h}. Now let us deduce \eqref{supremum_func_chi}. We have
\begin{eqnarray*}
\alpha |\ln(\alpha)|^p \Psi_{p,\alpha}( \alpha | \ln(\alpha)|) = \frac{|\ln(\alpha)|^p \times |\ln(\alpha |\ln(\alpha)|)|^{2-p}}{|\ln(\alpha)| + |\ln(\alpha |\ln(\alpha)|)|^2}  <  |\ln(\alpha)|^p \times |\ln(\alpha |\ln(\alpha)|)|^{-p}
\end{eqnarray*}
With the change of variable $\varrho = |\ln(\alpha)|$ (i.e. $\alpha = \exp{(- \varrho)}$), we have 
\begin{eqnarray*}
|\ln(\alpha)|^p \times |\ln(\alpha |\ln(\alpha)|)|^{-p} & =&   \frac{\varrho^p}{|\ln(\varrho \exp(-\varrho) )|^p}  \\
 & = & \frac{\varrho^p}{|-\varrho + \ln(\varrho)|^p} \\
 & = & \frac{\varrho^p}{(\varrho - \ln(\varrho))^p}  \to 1 \quad \text{as} \quad \varrho \to \infty. 
\end{eqnarray*}
This proves that 
$$
\alpha |\ln(\alpha)|^p \Psi_{p,\alpha}( \alpha | \ln(\alpha)|) = \bigO{1}
$$
and thus from \eqref{estimate root of h}, we deduce that \eqref{supremum_func_chi} holds.

\QED

\textbf{Proof of Proposition \ref{Prop order optimal conv rates under posteriori rule}.}
For simplicity of notation, let $\alpha := \alpha(\delta,y^\delta)$. In order to establish \eqref{eq tnnggggoifng}, we are going to bound the terms $||x^\dagger - x_\alpha||$ and $||x_\alpha - x_\alpha^\delta||$ separately. Let us start with the regularization error term. Given that $x^\dagger \in X_{f_p}(\rho)$, we have $x^\dagger =  f_p(T^*T) w$ and thus $x^\dagger - x_\alpha = r_\alpha(T^*T) x^\dagger = f_p(T^*T) r_\alpha(T^*T) w$. Hence by applying \cite[Proposition 1]{reg_exp_ipp} to $x^\dagger - x_\alpha$, we get
\begin{equation}
\label{eq xx}
|| x^\dagger - x_\alpha|| \leq || r_\alpha(T^*T)w|| \sqrt{\phi_p^{-1}\left(|| y -T x_\alpha||^2/\rho^2 \right)} \leq \rho \sqrt{\phi_p^{-1}\left(|| y -T x_\alpha||^2/\rho^2 \right)}.
\end{equation}
From \eqref{def_inv_func_phi_p} and \eqref{eq xx}, we deduce that 
\begin{equation}
\label{eq yy}
|| x^\dagger - x_\alpha|| \leq \rho f_p \left(|| y -T x_\alpha||^2/\rho^2  \right) (1 + o(1)).
\end{equation}
But 
\begin{eqnarray}
\label{eq zz}
|| y -T x_\alpha|| & \leq & || y^\delta -T x_\alpha^\delta|| + || y - T x_\alpha - (y^\delta - T x_\alpha^\delta ) || \nonumber \\ 
& \leq & \delta + \sqrt{\delta} + || r_\alpha(T^*T)(y - y^\delta)|| \nonumber\\
& \leq & 2\delta+\sqrt{\delta} \nonumber \\
& = & \sqrt{\delta}(2\sqrt{\delta}+1).
\end{eqnarray}
From \eqref{eq yy} and \eqref{eq zz}, we deduce that 
\begin{equation}
\label{eq tt}
|| x^\dagger - x_\alpha|| \leq \rho f_p \left( \delta (2\sqrt{\delta} +1)^2/\rho^2 \right)(1 + o(1)).
\end{equation}
Using \eqref{eq tt} and the fact that
\begin{equation}
\label{eq uu}
\frac{f_p \left( \delta (2\sqrt{\delta} +1)^2/\rho^2 \right)}{f_p(\delta)} = \left(\frac{-\ln(\delta)}{-\ln(\delta) - 2\ln(1+2\sqrt{\delta}) + 2\ln(\rho)} \right)^{p}  \to 1 \quad \text{as} \quad \delta \to 0,
\end{equation}
yields
\begin{equation}
\label{eq ii}
|| x^\dagger - x_\alpha|| = \bigO{f_p(\delta)} \quad \text{as}\quad \delta \to 0.
\end{equation}
Now let us estimate the propagated data noise term. Let $\bar{\alpha} = q \alpha$ with $q \in (1,2)$. From \eqref{morozov_like_par_sel_rule}, since $\bar{\alpha} > \alpha$, we get
\begin{equation}
\label{eq ff}
|| T x_{\bar{\alpha}}^\delta - y^\delta|| > \delta + \sqrt{\delta}.
\end{equation}
Therefore,
\begin{eqnarray}
\label{eq 1111}
|| T x_{\bar{\alpha}} - y|| & \geq & || T x_{\bar{\alpha}}^\delta - y^\delta|| - || T (x_{\bar{\alpha}}^\delta - x_{\bar{\alpha}}) - (y^\delta -y)|| \nonumber \\
& >& \delta + \sqrt{\delta} - || r_{\bar{\alpha}}(T^*T) (y^\delta - y)|| \nonumber \\
&\geq &\delta + \sqrt{\delta} - \delta \nonumber \\
&=& \sqrt{\delta}.
\end{eqnarray}
On the other hand, $|| T x_{\bar{\alpha}} - y||  = || T(x_{\bar{\alpha}} - x^\dagger)|| 
=|| (T^*T)^{1/2}(x_{\bar{\alpha}} - x^\dagger)|| =|| (T^*T)^{1/2} r_{\bar{\alpha}}(T^*T) x^\dagger|| $.
By applying \eqref{eq bound.} with $\varphi(t) = \sqrt{t}$ and $\epsilon = 1/8$, we get that there exists a constant $C$ such that
$|| (T^*T)^{1/2} r_{\bar{\alpha}}(T^*T) x^\dagger|| \leq C \bar{\alpha}^{3/8}$. This implies that
\begin{equation}
\label{eq noun}
|| T x_{\bar{\alpha}} - y|| \leq C \bar{\alpha}^{3/8}.
\end{equation}
From \eqref{eq 1111} and \eqref{eq noun}, we deduce that $\bar{\alpha}^{3/8} \geq \sqrt{\delta}/C$ which implies that $\bar{\alpha} \geq \bar{C} \delta^{4/3}$ with $\bar{C}= C^{-8/3}$.
From \eqref{bound fct prop err}, \eqref{bound prop error}, the above lower bound of $\bar{\alpha}$ and the fact that $\alpha > \bar{\alpha}/2$, we get that, there exists a positive constant $C'$ such that
\begin{equation}
\label{eq zzzt}
||x_\alpha - x_\alpha^\delta|| \leq C' \frac{\delta}{\sqrt{\alpha}} \leq C'\sqrt{2} \frac{\delta}{\sqrt{\bar{\alpha}}} \leq C'\sqrt{2/\bar{C}} \frac{\delta}{\sqrt{\delta^{4/3}}} = \delta^{1/3} C'\sqrt{2/\bar{C}}. 
\end{equation}
Given that $\delta^{1/3} = \bigO{f_p(\delta)}$ as $\delta \to 0$, we deduce that $||x_\alpha - x_\alpha^\delta|| = \bigO{f_p(\delta)}$ as $\delta \to 0$ which together with \eqref{eq ii} implies \eqref{eq tnnggggoifng}.
\QED
%\newpage

\end{document}